\documentclass[10pt]{article}
\usepackage{amssymb}
\usepackage{amsfonts}
\usepackage{color,graphicx}
\usepackage{latexsym}
\usepackage{amssymb}
\newtheorem{lem}{Lemma}[section]
\newtheorem{thm}[lem]{Theorem}

\newcommand{\proof}{\noindent {\it Proof.\hspace{4mm}}}
\newcommand{\qfd}{\ \framebox{}}

\def\newpic#1{%
\def\emline##1##2##3##4##5##6{%
\put(##1,##2){\special{em:point #1##3}}%
\put(##4,##5){\special{em:point #1##6}}%
\special{em:line #1##3,#1##6}}}
\newpic{}
\def\emline#1#2#3#4#5#6{%
\put(#1,#2){\special{em:moveto}}%
\put(#4,#5){\special{em:lineto}}}
\def\newpic#1{}
\newcommand\QQ{\hbox{I\kern-.53em\hbox{Q}}}
\newcommand{\Z}{\hbox{\bf Z}}

\title{On a $K_4$-UH self-dual 1-configuration $(102_4)_1$}
\author{Italo J. Dejter
\\ University of Puerto Rico \\ Rio Piedras, PR 00936-8377 \\ italo.dejter@gmail.com}
\date{}
\begin{document}
\maketitle

\begin{abstract}\noindent
Self-dual 1-con\-fi\-gu\-ra\-tions $(n_d)_1$ possess their Men\-ger
graph $\mathcal Y$ most $K_4$-separated among connected self-dual
con\-fi\-gu\-ra\-tions $(n_d)$. Such $\mathcal Y$ is most
sym\-me\-tric if $K_d$-ul\-tra\-ho\-mo\-ge\-neous. In this work,
such a $\mathcal Y$ is presented for $(n,d)=(102,4)$ and shown to
relate $n$ copies of the cuboctahedral graph $L(Q_3)$ to the $n$
copies of $K_d$; these are shown to share each copy of $K_3$ exactly
with two copies of $L(Q_3)$.
\end{abstract}

%\noindent{\bf Keywords:} Biggs-Smith graph, distance-transitive graph, ul\-tra\-ho\-mo\-ge\-neous graph, self-dual configuration, Menger graph

%\noindent{\bf 2000 Mathematics subject classification:} 05C62, 05B30, 05C20, 05C38 52C30, 68-04

\section{Introduction}

\noindent
%Let $1<c<m\in\Z$ and $1<d<n\in\Z$.
A {\it configuration} $R=(m_c,n_d)$ is an incidence structure of $m$
points and $n$ lines such that there are $c$ lines through each
point and $d$ points on each line \cite{Cox}. Thus, $cm = dn$. Let
$L=L(R)=L(m_c,n_d)$ be the bipartite graph with: {\bf(a)} $m$
``black'' vertices representing the points of $R$; {\bf(b)} $n$
``white'' vertices representing the lines of $R$; and {\bf(c)} an
edge between each two vertices representing a point and a line that
are incident in $R$. We call $L$ the {\it Levi graph} of $R$. If
$m=n$ and $c=d$, in which case $R$ is {\it symmetric}, then with
each configuration $R$ the dual configuration $\overline{R}$ may be
associated by reversing the roles of points and lines in $R$. Both
$R$ and $\overline{R}$ share the same Levi graph, but the
black-white coloring of their vertices is reversed. If $R$ is
isomorphic to its dual $\overline{R}$, then $R$ is {\it self-dual},
a corresponding isomorphism is called a {\it duality} and we denote
$R=(n_d)$. To any such configuration $(n_d)$ we can associate its
{\it Menger graph}, in which the points of $(n_d)$ are represented
by vertices, each two joined by an edge whenever the two
corresponding points are in a common line in $(n_d)$. If any two
points of $R$ are in at most $\lambda$ lines, then $R$ is a
$\lambda$-{\it configuration} $(n_d)_\lambda$ \cite{Gropp}. We note
from \cite{Cox} that the $4$-cube $Q_4$ may be considered as the
Levi graph of the M\"obius $(8_4)_2$ with "white" (resp. "black")
vertices being those of even (resp. odd) weight, (and so on for the
remaining Cox 2-configurations, in relation to the respective
$d$-cube $Q_d$).

\noindent Let $H$ be a connected regular graph. A graph $G$ is
$\mathcal C$-{\it ul\-tra\-ho\-mo\-ge\-neous} \cite{I}, or $\mathcal
C$-{\it UH}, if every isomorphism between two induced copies of
$H\in\mathcal C$ in $G$ extends to an auto\-mor\-phism of $G$. If
${\mathcal C}=\{H\}$ then $G$ is said to be $H$-{\it UH}. The
motivation of this paper is the study of connected Menger graphs
\cite{Cox} of self-dual 1-configurations $(n_d)_1$
\cite{handbook,Gropp} expressible as
$K_d$-ul\-tra\-ho\-mo\-ge\-neous graphs \cite{I}. The question of
for which values of $n$ such graphs exist is interesting because it
would yield the most symmetrical, connected, edge-disjoint unions of
$n$ copies of $K_d$ on $n$ vertices in which the roles of vertices
and copies of $K_d$ are interchangeable. For $d=4$, known values of
$n$ are: $n=13$, $21$ (see \cite{librofinal,gr,PS}) and $n=42$
(see\cite{D1}). While it is of interest to determine the spectrum
and multiplicities of the involved values of $n$, Theorem 4.1 below
contributes the value of $n=102$ via the Biggs-Smith association
scheme \cite{BCN}, later shown in Theorem 6.1 to control attachment
of 102 (cuboctahedral) copies of $L(Q_3)$ to the 102 (tetrahedral)
copies of $K_4$, these sharing each (triangular) copy of $K_3$ with
two copies of $L(Q_3)$ and guaranteeing in Theorem 7.1 that we have
the distance 3-graph of the Biggs-Smith graph $\mathcal S$
\cite{Bi,F} as the Menger graph $\mathcal Y$ of a self-dual
1-configuration $(102_4)_1$.

\noindent On the other hand, the Mobi\"us 2-configuration $(8_4)_2$
for example, and more generally the Cox 2-configurations
$((2^{d-1})_d)_2$ \cite{Cox}, have their Menger graphs with copies
of $K_4$ and $K_d$ respectively not edge-disjoint, even though these
are $K_4$- and $K_d$-ul\-tra\-ho\-mo\-ge\-neous graphs. Some
questions arising at this level are: Are variations of the latter
graphs as in \cite{PS} (5.3.7) $K_d$-ul\-tra\-ho\-mo\-ge\-neous?
Does there exist a relation between $K_d$-ul\-tra\-ho\-mo\-ge\-neous
Menger graphs and geometric configurations \cite{BP}? Do there exist
two different configurations with common $K_d$-ultrahomogeneous
Menger graph? Must $K_d$-ultrahomogeneous duality be involutory
\cite{GS,PS}?

\noindent A connected graph $G$ is an $\{H\}_n^d$-{\it graph} if it
is an edge-disjoint union of $n$ induced copies of $H$ with no other
copies of $H$ as subgraphs and each vertex incident to exactly $d$
copies of $H$, no two such copies sharing more than one vertex. If
$H=K_r$ is the complete graph of order $r$ ($0<r\in\Z$) then the
vertices and copies of $H$ in $G$ can be seen as the points and
lines of a 1-configuration $R_G$ with its points representing the
vertices of $G$ and its lines representing the copies of $H$ in $G$.
If $R_G$ is a self-dual 1-configuration, then it can be denoted
$(n_d)_1$ and $G$ can be recovered as the Menger graph of
$R_G=(n_d)_1$.

\noindent Let us illustrate these concepts with some examples.
Clearly, a connected graph $G$ is $m$-regular if and only if it is a
$\{K_2\}_{|E(G)|}^m$-graph. In this case, $G$ is arc-transitive if
and only if $G$ is $\{K_2\}$-UH. On the other hand:

{\bf(A)} for $1<r\in\Z$, the complete graph $K_r$ and its Cartesian
powers $K_r^2=K_r\,\square\,K_r, K_r^3=K_r^2\square
K_r,\ldots,K_r^s=K_r^{s-1}\,\square\,K_r,\ldots$ etc. are $K_r$-UH
$\{K_r\}_n^m$-graphs; their orders form a sequence
$r,r^2,r^3,\ldots,r^s,\ldots$ of integers corresponding to the
respective $K_r$-UH $\{K_r\}_1^1$-, $\{K_r\}_{2r}^2$-,
$\{K_r\}_{3r^2}^3$-, $\ldots$, $\{K_r\}_{sr^{s-1}}^r$-,
$\ldots$-graphs;

{\bf(B)} for $3\le r\in\Z$ the line graph $L(Q_r)$ of the $r$-cube
$Q_r$ is a $\{K_r,K_{2,2}\}$-UH
$\{K_r\}_{2^r}^2\{K_{2,2}\}_{r(r-1)2^{r-3}}^{r-1}$-graph. A similar
argument yields a $K_r$-UH $\{K_r\}_n^m$-graph out of any other
regular-polytopal graph via its line graph.

\noindent There is  only one case in {\bf(A)}-{\bf(B)} that is
Menger graph of a self-dual configuration, namely $K_2^2$ (duality
sending for example the points $00,10,11,01$ resp. onto the lines
$x0,0x,x1,1x$, where $0\le x\le 1$), even though all graphs $K_r^r$
have equal numbers of vertices and of copies of $K_r$ so they are
Menger graphs of {\it symmetric} configurations (but not self-dual).
If $r=4$, then the orders of the $K_d$-UH $\{K_d\}_n^m$-graphs in
{\bf(A)}-{\bf(B)} are divisible by 4. Beside ours ($n=132$), a case
of even order indivisible by 4 is the one mentioned above on $n=42$
vertices \cite{D1}. Its construction was based on the ordered
pencils of the Fano plane. Extensions of that construction of
\cite{D1}, based on ordered pencils of binary projective spaces, are
introduced in \cite{Dinf}, which provides $K_4$-UH
$\{K_4\}_n^m$-graphs whose even orders are indivisible by 4, the
smallest of which being 210. However, the latter graphs are not
Menger graphs of self-dual configurations. A configuration $(n_d)_1$
is said to be $K_d$-{\it UH} if its Menger graph is. Are there any
UH-$K_4$ self-dual configurations $(n_4)_1$ with even $n<42$? Or
$42<n<102$?

\noindent In Section 4, the claimed Menger graph $\mathcal Y$ is
constructed by means of the distance-3 graphs of the 9-cycles of the
Biggs-Smith graph $\mathcal S$. Theorem 4.1 proves our claim about
$\mathcal Y$ as an application of a transformation of
distance-transitive graphs into $\mathcal C$-UH graphs that took in
\cite{DCox} from the Coxeter graph of order 28 onto the Klein graph
of order 56. A similar application allowed in \cite{D-PD} to
confront, as digraphs, the Pappus graph of order 18 to the Desargues
graph of order 20. These applications as well as \cite{orestes} use
the following definitions. Given a family $\mathcal C$ of digraphs,
a digraph $G$ is said to be $\mathcal C$-{\it UH} if every
isomorphism between two induced members of $\mathcal C$ in $G$
extends to an auto\-mor\-phism of $G$. If ${\mathcal C}=\{H\}$ then
$G$ is said to be $H$-{\it UH}. By removing the suffix ``di'' here,
the definition of $\mathcal C$-UH graph is recovered. A presentation
of $\mathcal S$ is given in Section 2 by means of Biggs-Hoare
sextets mod 17 \cite{BH} which provide a convenient notation to
present $\mathcal Y$ in Section 3 in preparation for Section 4.
%Sections 5, 6 and 7 employ the Biggs-Smith association scheme in order to prove the remaining claims on $\mathcal Y$.

\noindent We set one more definition, to be used from Section 2 on.
If $M$ is a subgraph of $H$ and if $G$ is both $M$-UH, and $H$-UH,
then $G$ is an $\{H\}_{M}$-{\it UH graph} if, for each induced copy
$H_0$ of $H$ in $G$ containing an induced copy $M_0$ of $M$, there
exists exactly one induced copy $H_1\neq H_0$ of $H$ in $G$ such
that $V(H_0)\cap V(H_1)=V(M_0)$ and $E(H_0)\cap E(H_1)=E(M_0)$.

\section{The Biggs-Smith graph}

%\newpage

\begin{figure}[htp]
\vspace*{-2mm}
\hspace*{15mm}
  \includegraphics[scale=0.20]{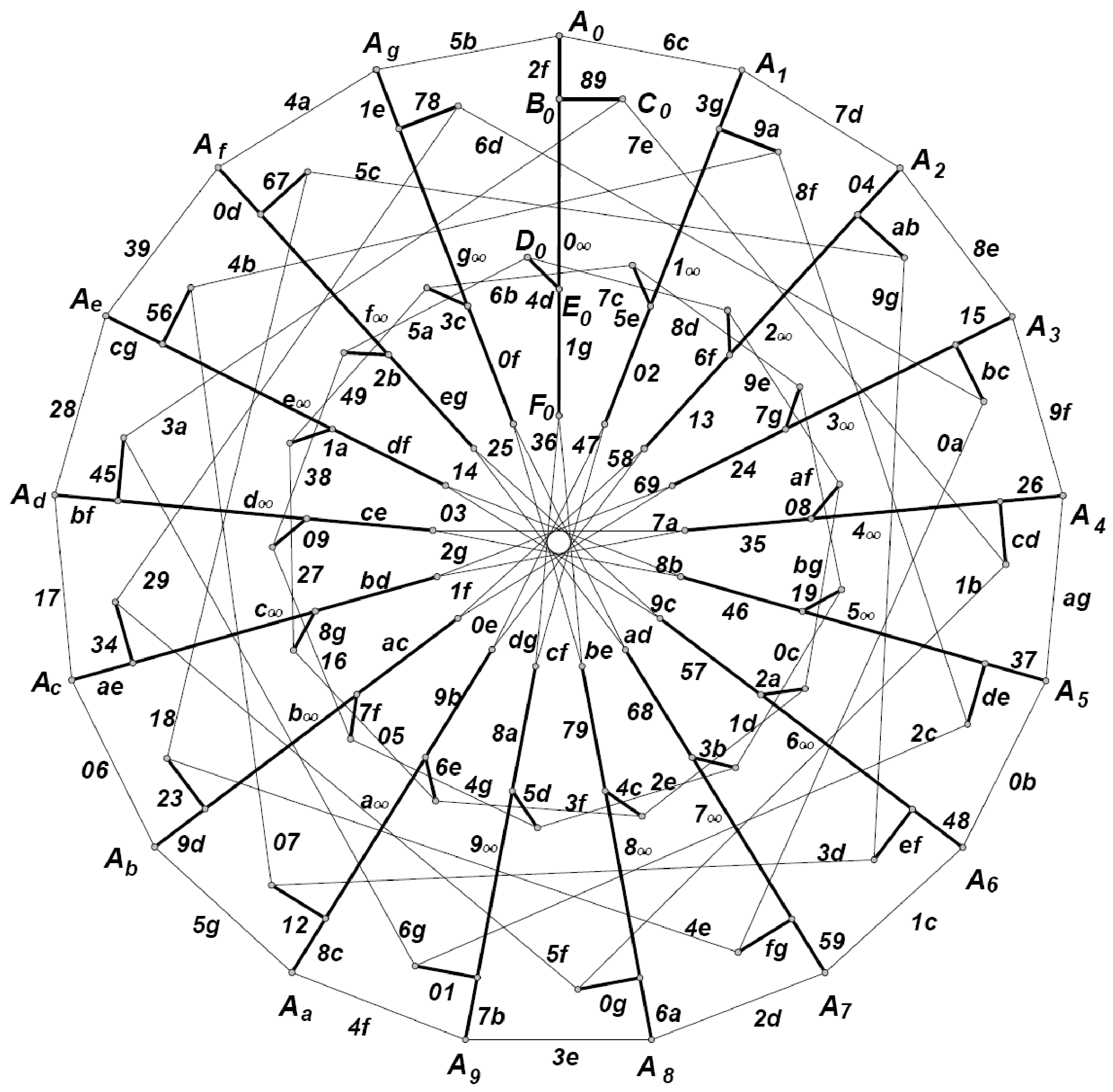}\\%0.43dancers
%\vspace*{-8mm}
\caption{Representation of $\mathcal S$ via sextets
and thick subtrees $T_i^\infty$}
\end{figure}

\noindent The Biggs-Smith graph $\mathcal S$ has order $n=102$,
diameter $d=7$, girth $g=9$ and automorphism group ${\mathcal
A}=PSL(2,17)$ \cite{BCN}. By letting $k$ be the largest integer $s$
such that $\mathcal S$ is $s$-arc transitive, it is seen that $k=4$.
In addition, the number $\eta$ of $9$-cycles of $\mathcal S$ is
$\eta=136$. Taking into account the definition in the last paragraph
of Section 1 and by denoting a 3-path by $P_4$ and a 9-cycle by
$\gamma_9$, the following particular case of Theorem 3 of
\cite{orestes} holds (which cannot be refined to a result of
$\{\vec{\gamma_9}\}_{\vec{P_4}}$-UH digraphs; see (4) below):
\begin{equation}{\mathcal S}\mbox{ is }\{\gamma_9\}_{P_4}\mbox{-{\rm
UH}}.\end{equation} Properties of $\mathcal S$ we need are presented
via {\it sextets} \cite{BH}, where heptadecimal notation is used to
denote elements of $GF(17)$ (for example $g=16=-1$ and $d=13=-4$) so
that $\mathcal S$ can be characterized as a connected graph whose
vertex set $V({\mathcal S})$ comprises 102 sextets mod 17, namely
102 unordered triples $$\{a_0b_0,a_1b_1,a_2b_2\}$$ composed by
unordered pairs $a_ib_i$ of points $a_i,b_i$ of the projective line
$PG(1,17)=GF(17)\cup\{\infty\}$ satisfying

$$(a_i-a_j)(b_i-b_j)(a_i-b_j)^{-1}(b_i-a_j)^{-1}=-1,$$ if $a_i\ne\infty$ and satisfying
$$(b_i-b_j)(b_i-a_j)^{-1}=-1,$$ if $a_i=\infty$, whenever $i\ne j$ in
$\{0,1,2\}$, including the vertices
\begin{eqnarray}\begin{array}{ccc}A_0=\{2f,5b,6c\},& B_0=\{0\infty,2f,89\},&
C_0=\{3a,7e,89\},\\
D_0=\{5a,7c,4d\},& E_0=\{0\infty,1g,4d\},&
F_0=\{1g,36,be\}.\end{array}\end{eqnarray} Any two of the resulting
102 vertices are adjacent in $\mathcal S$ whenever they share one
such pair $a_ib_i$, in which case the resulting edge is labeled
$a_ib_i$. It is shown in \cite{BH} that this $\mathcal S$ is unique
and that the {\it edge labels} $a_ib_i$ are pairwise distinct, so
they determine an edge labeling of $\mathcal S$ represented in
Figure 1 with the following notation. The six vertices in (2) are
those of a subtree $T_0^\infty$ (of $\mathcal S$) which is the
edge-disjoint union of the paths
$$(A_0,2f,B_0,89,C_0),(D_0,4d,E_0,1g,F_0)\mbox{ and }(B_0,0\infty,E_0)$$
of lengths 3, 3 and 2, respectively. By adding to all elements of
$GF(17)$ in $T_0^\infty$ a constant $i\in GF(17)$, a similar tree
$T_i^\infty$ is obtained. The trees $T_0^\infty,\ldots,T_g^\infty$,
represented in Figure 1 via dark traces, are pairwise disjoint and
cover $V({\mathcal S})$. The complement of their union in $\mathcal
S$ is formed by $4$ 17-cycles
$$\begin{array}{ll}A=(A_0$,6c,$A_1, \ldots, A_g,5b),&
D=(D_0,7c,D_2, \ldots, D_f,5a),\\ C=(C_0,7e,C_4, \ldots, C_d,3a),&
F=(F_0,be,F_8, \ldots, F_9,36).\end{array}$$ Each of these cycles
$y=A,D,C,F$ has vertices $y_r$ with $r\in GF(17)$ advancing in 1, 2,
4, 8 units mod 17 stepwise from left to right, respectively.

\noindent Employed in \cite{orestes} in proving (1) above, there is
a set ${\mathcal C}_9$ of 136 directed 9-cycles of $\mathcal S$, of
which a generating subset
$$\{\Pi^0=(\Pi_0^0\Pi_1^0\ldots\Pi_8^0);\Pi=S,T,\ldots,Z\}$$ (written
without commas and accompanied to the right by auxiliary
permutations, as explained below) is as follows:
\begin{equation}\begin{array}{ll}
^{S^0=(B_2A_2A_1A_0A_gA_fB_fC_fC_2)}
_{T^0=(E_dD_dD_fD_0D_2D_4E_4F_4F_d)}& ^{s^0=(07cb4d65a)(\infty
8g2e3f19)}
_{t^0=(03ac9857e)(\infty12d6b4fg)}\\
^{U^0=(B_9C_9C_dC_0C_4C_8B_8A_8A_9)}
_{V^0=(E_gF_gF_8F_0F_9F_1E_1D_1D_g)}& ^{u^0=(06371gaeb)(\infty
249c58df)}
_{v^0=(05b3f2e6c)(\infty d9ga7184)}\\
^{W^0=(B_9E_9F_9F_0F_8E_8B_8A_8A_9)}
_{X^0=(E_gB_gA_gA_0A_1B_1E_1D_1D_g)}& ^{w^0=(\infty
a3b986e7)(0df15cg24)}
_{x^0=(\infty ebcg1563)(084f7a2d9)}\\
^{Y^0=(B_2E_2D_2D_0D_fE_fB_fC_fC_2)}
_{Z^0=(E_dB_dC_dC_0C_4B_4E_4F_4F_d)}& ^{y^0=(\infty
6ca2f75b)(01943ed8g)}
_{z^0=(\infty 5aed437c)(0fg9b6812)}\\
\end{array}\end{equation}
where the permutation
$\pi^0=(\pi_0^0\pi_1^0\ldots\pi_8^0)(\xi_0^0\xi_1^0\ldots\xi_8^0)$
of $PG(1,17)$ to the right of each $\Pi^0$ is such that: {\bf(i)}
the pair $\pi_i^0\pi_{i+4}^0$ labels the edge $\Pi_i^0\Pi_{i+1}^0$;
{\bf(ii)} the pair $\xi_i^0\xi_{i+3}^0$ labels the only edge
incident to $\Pi_i^0$ outside $\Pi^0$, where $i=0,\ldots,8$ and
index addition is taken modulo 9. ${\mathcal C}_9$ also contains the
directed cycles $\Pi^r$ with accompanying permutations $\pi^r$
obtained from $\Pi^0$ and $\pi^0$ by uniformly adding $r\in\Z_{17}$
mod 17 to all subscripts and superscripts. Observe that: {\bf(iii)}
passing from $s^0$ to $t^0$ to $u^0$ to $v^0$ and again to $s^0$,
(resp. from $w^0$ to $x^0$ to $y^0$ to $z^0$ and again to $w^0$)
amounts to multiplying uniformly and successively the participating
entries of the permutations $\pi^0$ by either $2$ or $-2$ mod 17;
and {\bf(iv)} $S^0,\ldots,Z^0$ are invariant with respect to their
change-of-sign involutions mod 17, with corresponding involutions on
$s^0,\ldots,z^0$ around the initial entries of their two composing
cycles, which are either 0 and $\infty$, or $\infty$ and 0.

%\section{Oriented 9-cycles and distance-3 digraphs}
\section{Distance-3 digraphs of oriented 9-cycles}

\noindent A $k$-{\it arc} in a (di)graph is a sequence of vertices
$v_0v_1\ldots v_k$ (written without parentheses or commas), where
consecutive vertices are adjacent and $v_{i-1}\ne v_{i+1}$, for
$0<i<k$ \cite{GR}. A $k$-arc can be interpreted as a directed walk
of length $k$ in which consecutive edges are distinct \cite{GY}.
Thus, an arc in a (di)graph $\Gamma$ is a 1-arc of $\Gamma$.  The
form in which the directed 9-cycles $\Pi^{r}$ in Section 2 share
3-arcs, either oppositely oriented or not, to be used in Figure 3
below, can be encoded as in the following table that for each
$\Pi^0$ presents details (explained below) of the 9-cycles
$\Xi_r\ne\Pi^0$ in ${\mathcal C}_9$ that intersect $\Pi^0$ either in
the succeeding 3-arcs $\Pi^0_i\Pi^0_{i+1}\Pi^0_{i+2}\Pi^0_{i+3}$ or
in their respective reversed arcs, for $i=0,\ldots,8$, with sums
involving $i$ taken mod 9:
\begin{equation}\begin{array}{c}
^{S^0\;:(\mbox{-}X^1_2,  \,\,\,S^1_2,
\,\,\,S^g_1,\mbox{-}X^g_1,\mbox{-}U^7_5,\,\,\,U^6_8,\,\,\,Y^0_6,
\,\,\,U^b_4,\mbox{-}U^a_7);} _{T^0\;:(\mbox{-}Y^f_2,\,\,\,T^f_2,
\,\,\,T^2_1,\mbox{-}Y^2_1,\mbox{-}V^3_5,\,\,\,V^5_8,\,\,\,Z^0_6,
\,\,\,V^c_4,\mbox{-}V^e_7);}\\
^{U^0\;:(  \,\,\,Z^d_1,\,\,\,U^d_2,\,\,\,U^4_1,\,\,\,Z^4_2,
\,\,\,S^6_7,\mbox{-}S^a_4,\,\,\,W^0_6,\mbox{-}S^7_8, \,\,\,S^b_5);}
_{V^0\;:(\mbox{-}W^8_2,\,V^8_2,
\,\,\,V^9_1,\mbox{-}W^9_1,\,T^5_7,\mbox{-}T^e_4,\,\,\,X^0_6,
\mbox{-}T^3_8,  \,\,\,T^c_5);}\\
^{W^0:(\mbox{-}Z^d_7,\mbox{-}V^8_3,\mbox{-}V^9_0,\mbox{-}Z^4_5,\mbox{-}W^g_8,
\,X^9_0,\,\,\,U^0_6,    \,\,\,X^8_3,\mbox{-}W^1_4);}
_{X^0\,:(\,W^8_5,\mbox{-}S^1_3,\mbox{-}S^g_0,\,\,\,W^9_7,
\mbox{-}X^2_8,\,\,\,Y^g_0,\,\,\,V^0_6,\,\,\,Y^1_3,\mbox{-}X^f_4);}\\
^{Y^0\,:(
\,\,\,X^1_5,\mbox{-}T^f_3,\mbox{-}T^2_0,\,\,\,X^g_7,\mbox{-}Y^d_8,
\,\,\,Z^2_0,\,\,\,S^0_6,\,\,\,\,\,Z^f_3,\mbox{-}Y^4_4);} _{Z^0\,:(
\,\,\,Y^f_5,  \,\,\,U^4_0,  \,\,\,U^d_3,
\,\,\,Y^2_7,\mbox{-}Z^8_8,\mbox{-}W^d_3,\,\,\,T^0_6,
\mbox{-}W^4_0,\mbox{-}Z^9_4).}
\end{array}\end{equation}
Each such $\Xi^r$ has: either {\bf(I)} a preceding minus sign, if
the corresponding 3-arcs in $\Pi^0$ and $\Xi^r$ are oppositely
oriented, or {\bf(II)} no preceding sign, otherwise. Each shown
$-\Xi^r_j$ (resp. $\Xi^r_j$) has a subscript $j$ indicating the
equality of initial vertices $\Xi^r_j=\Pi^0_{i+3}$ (resp.
$\Xi^{\,r}_j=\Pi^0_i$) of those 3-arcs, for $i=0,\ldots,8$.

\noindent Given a (di)graph $\Gamma$ and a positive integer $k\le$
diameter$(\Gamma)$, the {\it distance}-$k$ ({\it di}){\it graph}
$\Gamma_k$ of $\Gamma$, with vertex set $V(\Gamma_k)=V(\Gamma)$, is
such that from every $u\in V(\Gamma_k)$ an arc of $\Gamma_k$ departs
to a vertex $v\ne u$ whenever there is a shortest $k$-arc of length
$k$ in $\Gamma$ from $u$ to $v$. Let $({\mathcal C}_9)_3$ be the
family of distance-3 digraphs of directed $9$-cycles in ${\mathcal
C}_9$. On a representation of an arc $e=w_0w_1$ of a member
$(\zeta_9)_3$ of $({\mathcal C}_9)_3$, we label its {\it tail}, or
initial vertex, $w_0$, its {\it initial flag} $\{w_0,e\}$, its {\it
terminal flag} $\{e,w_1\}$ and its {\it head}, or terminal vertex,
$w_1$, respectively by the names of the vertices $v_0,v_1,v_2,v_3$
of the 3-arc $v_0v_1v_2v_3$ in $\zeta_9$ for which $w_0w_1$ stands
in $(\zeta_9)_3$. For example, if
$\zeta_9=U^9=(B_1C_1C_5C_9C_dC_0B_0A_0A_1)$, so that $(\zeta_9)_3=$
$(U^9)_3=$ $(B_1C_9B_0)(C_1C_d$ $A_0)(C_5C_0A_1)$, then the initial
flag of the arc $B_1C_9$ in $(\zeta_9)_3=(U^9)_3$ is labeled by
$C_1$, the terminal flag by $C_5$, while $B_1$ and $C_9$ are labeled
exactly by $B_1$ and $C_9$, respectively. We get the labels over
$(\zeta_9)_3=(U^0)_3$ shown in Figure 2.

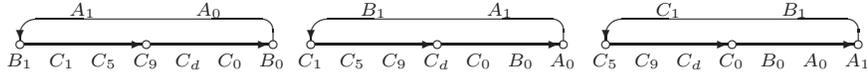
\begin{figure}[htp] %\hspace*{18mm}
\unitlength=0.56mm \special{em:linewidth 0.4pt}
\linethickness{0.4pt}
\begin{picture}(201.00,13.00)
\put(1.00,5.00){\circle{2.00}} \put(31.00,5.00){\circle{2.00}}
\put(1.00,1.00){\makebox(0,0)[cc]{$_{B_1}$}}
\put(21.00,1.00){\makebox(0,0)[cc]{$_{C_5}$}}
\put(31.00,1.00){\makebox(0,0)[cc]{$_{C_9}$}}
\put(61.00,5.00){\circle{2.00}}
\put(11.00,1.00){\makebox(0,0)[cc]{$_{C_1}$}}
\put(51.00,1.00){\makebox(0,0)[cc]{$_{C_0}$}}
\put(61.00,1.00){\makebox(0,0)[cc]{$_{B_0}$}}
\put(41.00,1.00){\makebox(0,0)[cc]{$_{C_d}$}}
\put(16.00,13.00){\makebox(0,0)[cc]{$_{A_1}$}}
\put(16.00,6.00){\oval(30.00,10.00)[lt]}
\put(46.00,13.00){\makebox(0,0)[cc]{$_{A_0}$}}
\put(70.00,5.00){\circle{2.00}} \put(100.00,5.00){\circle{2.00}}
\put(70.00,1.00){\makebox(0,0)[cc]{$_{C_1}$}}
\put(90.00,1.00){\makebox(0,0)[cc]{$_{C_9}$}}
\put(100.00,1.00){\makebox(0,0)[cc]{$_{C_d}$}}
\put(130.00,5.00){\circle{2.00}}
\put(80.00,1.00){\makebox(0,0)[cc]{$_{C_5}$}}
\put(120.00,1.00){\makebox(0,0)[cc]{$_{B_0}$}}
\put(130.00,1.00){\makebox(0,0)[cc]{$_{A_0}$}}
\put(110.00,1.00){\makebox(0,0)[cc]{$_{C_0}$}}
\put(85.00,13.00){\makebox(0,0)[cc]{$_{B_1}$}}
\put(85.00,6.00){\oval(30.00,10.00)[lt]}
\put(115.00,13.00){\makebox(0,0)[cc]{$_{A_1}$}}
\put(140.00,5.00){\circle{2.00}} \put(170.00,5.00){\circle{2.00}}
\put(140.00,1.00){\makebox(0,0)[cc]{$_{C_5}$}}
\put(160.00,1.00){\makebox(0,0)[cc]{$_{C_d}$}}
\put(170.00,1.00){\makebox(0,0)[cc]{$_{C_0}$}}
\put(200.00,5.00){\circle{2.00}}
\put(150.00,1.00){\makebox(0,0)[cc]{$_{C_9}$}}
\put(190.00,1.00){\makebox(0,0)[cc]{$_{A_0}$}}
\put(200.00,1.00){\makebox(0,0)[cc]{$_{A_1}$}}
\put(180.00,1.00){\makebox(0,0)[cc]{$_{B_0}$}}
\put(155.00,13.00){\makebox(0,0)[cc]{$_{C_1}$}}
\put(155.00,6.00){\oval(30.00,10.00)[lt]}
\put(185.00,13.00){\makebox(0,0)[cc]{$_{B_1}$}}
\put(46.50,8.00){\oval(29.00,6.00)[rt]}
\emline{61.00}{8.00}{1}{61.00}{6.00}{2}
\put(115.50,8.00){\oval(29.00,6.00)[rt]}
\emline{130.00}{8.00}{3}{130.00}{6.00}{4}
\put(185.50,8.00){\oval(29.00,6.00)[rt]}
\emline{200.00}{8.00}{5}{200.00}{6.00}{6}
\emline{155.00}{11.00}{7}{186.00}{11.00}{8}
\emline{85.00}{11.00}{9}{116.00}{11.00}{10}
\emline{16.00}{11.00}{11}{47.00}{11.00}{12}
\put(2.00,5.00){\vector(1,0){28.00}}
\put(32.00,5.00){\vector(1,0){28.00}}
\put(71.00,5.00){\vector(1,0){28.00}}
\put(101.00,5.00){\vector(1,0){28.00}}
\put(1.00,7.00){\vector(0,-1){1.00}}
\put(70.00,7.00){\vector(0,-1){1.00}}
\put(140.00,7.00){\vector(0,-1){1.00}}
\put(141.00,5.00){\vector(1,0){28.00}}
\put(171.00,5.00){\vector(1,0){28.00}}
\end{picture}
\caption{Labels of vertices and flags of $(\zeta_9)_3=(U^9)_3$}
\end{figure}

\section{$K_4$-UH self-dual 1-configuration $(102_4)_1$}

\noindent We are to fasten pairs of arcs of the digraphs
$(\zeta_9)_3$ defined in Section 3 in such a way that a graph
$\mathcal Y$ with the properties claimed in Section 1 is produced. A
sequence of operations ${\mathcal S}\rightarrow{\mathcal
C}_9\rightarrow({\mathcal C}_9)_3\rightarrow{\mathcal Y}$ (compare
with \cite{DCox}) is performed in order to transform $\mathcal S$
into the claimed $\mathcal Y$. Each distance-3 digraph $(\zeta_9)_3$
of a 9-cycle $\zeta_9$ in the collection ${\mathcal C}_9$ generated
via (3) is formed by 3 disjoint directed triangles. It yields a
total of $3\times 136$ directed triangles so ${\mathcal C}_9$
determines a family of 408 directed triangles in the claimed
$\mathcal Y$ with each edge shared by exactly two such directed
triangles in arcs that are either oppositely or identically
oriented. It amounts to 102 copies of $K_4$; these can be subdivided
into 6 subfamilies $\{\Sigma^i\}$ of 17 copies each, say with
$\Sigma\in\{A,B,C,D,E,F\}$ and $i\in\{0,1,\ldots,16=g\}=\Z_{17}$.
The vertex sets $V(\Sigma^i)$, each followed by the set
$\Lambda(\Sigma_i)$ of copies of $K_4$ containing the corresponding
vertex $\Sigma_i$ can be taken as follows, showing $\Z_2$-symmetry
produced by change of sign mod 17:
\begin{equation}\begin{array}{lllll}
^{V(A^i)=\{C_i,}
_{V(B^i)=\{D_{i+3},}&
^{D_i,}
_{D_{i-3},}&
^{E_{i+4},E_{i-4}\};\Lambda(A_i)=\{C^i,}
_{F_{i+5},F_{i-5}\};\Lambda(B_i)=\{D^{i+2},}&
^{D^i,}
_{D^{i-2},}&
^{E^{i+7},E^{i-7}\};}
_{F^{i+8},F^{i-8}\};}\\
^{V(C^i)=\{A_i,}_
{V(D^i)=\{A_i,}&
^{F_i,}_{D_i,}&
^{E_{i+1},E_{i-1}\};\Lambda(C_i)=\{A^i,}
_{B_{i+2},B_{i-2}\};\Lambda(D_i)=\{A^i,}&
^{F^i,}
_{D^i,}&
^{E^{i+6},E^{i-6}\};}
_{B^{i+3},B^{i-3}\};}\\
^{V(E^i)=\{C_{i+6},}_{V(F^i)=\{C_i,}&^{C_{i-6},}_{F_i,}&
^{A_{i+7},A_{i-7}\};\Lambda(E_i)=\{C^{i+1},}_{B_{i+8},B_{i-8}\};\Lambda(F_i)=\{C^i,}&^{C^{i-1},}_{F^i,}&
^{A^{i+4},A^{i-4}\};}_{B^{i+5},B^{i-5}\};}
\end{array}\end{equation}
\noindent where $i$ varies in $\Z_{17}$. This reveals a duality
$\phi$ from the 102 vertices of $\mathcal S$ onto the 102 copies of
$K_4$ in $\mathcal S$. In fact, these copies of $K_4$ are the
vertices of a graph $\phi({\mathcal S})={\mathcal
S}^*\equiv{\mathcal S}$ determined by
\begin{equation}\begin{array}{lll}^{\phi(A_i)=A^{3i}=A_i^*,}_{\phi(D_i)=D^{5i}=D_i^*,}&
^{\phi(B_i)=B^{-7i}=B_i^*,}_{\phi(E_i)=E^{6i}\;\;=E_i^*,}&^{\phi(C_i)=C^{3i}=C_i^*,}_{
\phi(F_i)=F^{5i}=F_i^*,}\end{array}\end{equation} ($i\in\Z_{17}$),
with a structure similar to that of the vertices $A_i,\ldots,F_i$ of
$\mathcal S$, the copies of $K_4$ in ${\mathcal S}^*$
precisely being $\Sigma_i=A_i,\ldots,F_i$ and corresponding vertex
sets $\Lambda(\Sigma_i)$ as specified above. Moreover,
$\phi:{\mathcal S}\rightarrow{\mathcal S}^*$ is a graph isomorphism,
with the adjacency of ${\mathcal S}^*$ equivalent to that of
$\mathcal S$.

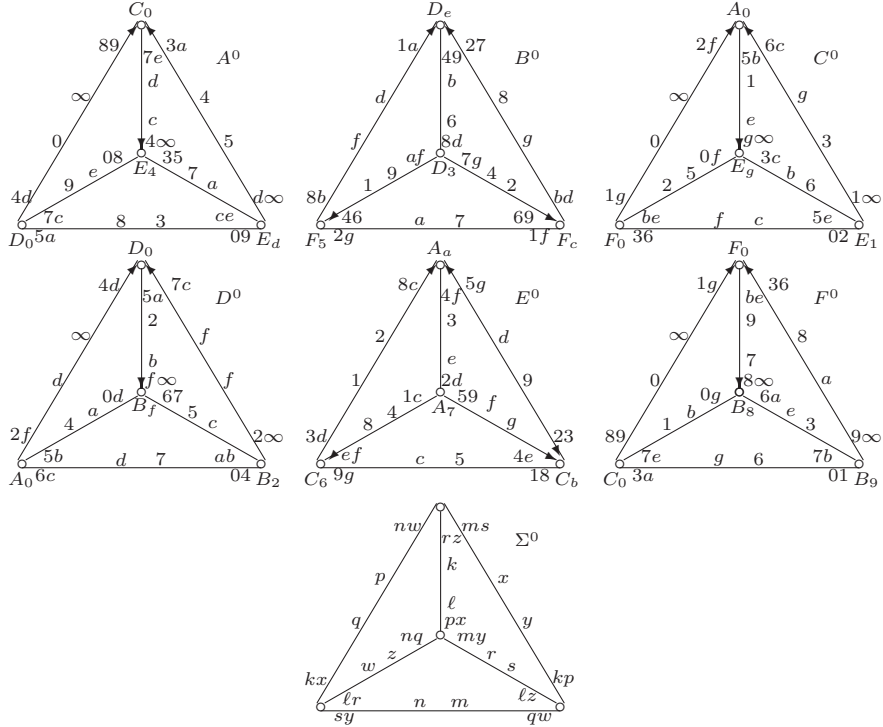
\begin{figure}[htp]
%\vspace*{-13mm}
%\input{4u3right.pic}
\unitlength=0.53mm \special{em:linewidth 0.4pt}
\linethickness{0.4pt}
\begin{picture}(213.00,179.00)
\put(151.00,125.00){\circle{2.00}}
\put(211.00,125.00){\circle{2.00}}
\put(181.00,175.00){\circle{2.00}}
\put(181.00,143.00){\circle{2.00}}
\emline{151.00}{124.00}{1}{211.00}{124.00}{2}
\emline{180.00}{176.00}{3}{150.00}{126.00}{4}
\emline{152.00}{126.00}{5}{180.00}{142.00}{6}
\put(151.00,65.00){\circle{2.00}} \put(211.00,65.00){\circle{2.00}}
\put(181.00,115.00){\circle{2.00}} \put(181.00,83.00){\circle{2.00}}
\emline{151.00}{64.00}{7}{211.00}{64.00}{8}
\emline{212.00}{66.00}{9}{182.00}{116.00}{10}
\emline{182.00}{82.00}{11}{210.00}{66.00}{12}
\put(76.00,125.00){\circle{2.00}} \put(136.00,125.00){\circle{2.00}}
\put(106.00,175.00){\circle{2.00}}
\put(106.00,143.00){\circle{2.00}}
\emline{137.00}{126.00}{13}{107.00}{176.00}{14}
\emline{105.00}{176.00}{15}{75.00}{126.00}{16}
\emline{77.00}{126.00}{17}{105.00}{142.00}{18}
\emline{107.00}{142.00}{19}{135.00}{126.00}{20}
\put(76.00,65.00){\circle{2.00}} \put(136.00,65.00){\circle{2.00}}
\put(106.00,115.00){\circle{2.00}} \put(106.00,83.00){\circle{2.00}}
\emline{137.00}{66.00}{21}{107.00}{116.00}{22}
\emline{105.00}{116.00}{23}{75.00}{66.00}{24}
\emline{77.00}{66.00}{25}{105.00}{82.00}{26}
\emline{107.00}{82.00}{27}{135.00}{66.00}{28}
\put(1.00,125.00){\circle{2.00}} \put(61.00,125.00){\circle{2.00}}
\put(31.00,175.00){\circle{2.00}} \put(31.00,143.00){\circle{2.00}}
\emline{1.00}{124.00}{29}{61.00}{124.00}{30}
\emline{30.00}{176.00}{31}{0.00}{126.00}{32}
\emline{2.00}{126.00}{33}{30.00}{142.00}{34}
\put(1.00,65.00){\circle{2.00}} \put(61.00,65.00){\circle{2.00}}
\put(31.00,115.00){\circle{2.00}} \put(31.00,83.00){\circle{2.00}}
\emline{1.00}{64.00}{35}{61.00}{64.00}{36}
\emline{62.00}{66.00}{37}{32.00}{116.00}{38}
\emline{32.00}{82.00}{39}{60.00}{66.00}{40}
\put(7.00,62.00){\makebox(0,0)[cc]{$_{6c}$}}
\put(9.00,67.00){\makebox(0,0)[cc]{$_{5b}$}}
\put(1.00,72.00){\makebox(0,0)[cc]{$_{2f}$}}
\put(7.00,122.00){\makebox(0,0)[cc]{$_{5a}$}}
\put(9.00,127.00){\makebox(0,0)[cc]{$_{7c}$}}
\put(1.00,132.00){\makebox(0,0)[cc]{$_{4d}$}}
\put(23.00,170.00){\makebox(0,0)[cc]{$_{89}$}}
\put(34.00,167.00){\makebox(0,0)[cc]{$_{7e}$}}
\put(40.00,170.00){\makebox(0,0)[cc]{$_{3a}$}}
\put(52.00,127.00){\makebox(0,0)[cc]{$_{ce}$}}
\put(39.00,142.00){\makebox(0,0)[cc]{$_{35}$}}
\put(36.00,146.00){\makebox(0,0)[cc]{$_{4\infty}$}}
\put(24.00,142.00){\makebox(0,0)[cc]{$_{08}$}}
\put(56.00,122.00){\makebox(0,0)[cc]{$_{09}$}}
\put(63.00,132.00){\makebox(0,0)[cc]{$_{d\infty}$}}
\put(10.00,146.00){\makebox(0,0)[cc]{$_0$}}
\put(16.00,157.00){\makebox(0,0)[cc]{$_\infty$}}
\put(34.00,161.00){\makebox(0,0)[cc]{$_d$}}
\put(34.00,151.00){\makebox(0,0)[cc]{$_c$}}
\put(47.00,157.00){\makebox(0,0)[cc]{$_4$}}
\put(53.00,146.00){\makebox(0,0)[cc]{$_5$}}
\put(13.00,135.00){\makebox(0,0)[cc]{$_9$}}
\put(19.00,138.00){\makebox(0,0)[cc]{$_e$}}
\put(49.00,135.00){\makebox(0,0)[cc]{$_a$}}
\put(44.00,138.00){\makebox(0,0)[cc]{$_7$}}
\put(36.00,126.00){\makebox(0,0)[cc]{$_3$}}
\put(26.00,126.00){\makebox(0,0)[cc]{$_8$}}
\put(31.00,179.00){\makebox(0,0)[cc]{$_{C_0}$}}
\put(63.00,121.00){\makebox(0,0)[cc]{$_{E_d}$}}
\put(1.00,121.00){\makebox(0,0)[cc]{$_{D_0}$}}
\put(32.00,139.00){\makebox(0,0)[cc]{$_{E_4}$}}
\put(23.00,110.00){\makebox(0,0)[cc]{$_{4d}$}}
\put(34.00,107.00){\makebox(0,0)[cc]{$_{5a}$}}
\put(41.00,110.00){\makebox(0,0)[cc]{$_{7c}$}}
\put(39.00,82.00){\makebox(0,0)[cc]{$_{67}$}}
\put(36.00,86.00){\makebox(0,0)[cc]{$_{f\infty}$}}
\put(24.00,82.00){\makebox(0,0)[cc]{$_{0d}$}}
\put(56.00,62.00){\makebox(0,0)[cc]{$_{04}$}}
\put(63.00,72.00){\makebox(0,0)[cc]{$_{2\infty}$}}
\put(10.00,86.00){\makebox(0,0)[cc]{$_d$}}
\put(16.00,97.00){\makebox(0,0)[cc]{$_\infty$}}
\put(34.00,101.00){\makebox(0,0)[cc]{$_2$}}
\put(34.00,91.00){\makebox(0,0)[cc]{$_b$}}
\put(47.00,97.00){\makebox(0,0)[cc]{$_f$}}
\put(53.00,86.00){\makebox(0,0)[cc]{$_f$}}
\put(13.00,75.00){\makebox(0,0)[cc]{$_4$}}
\put(19.00,78.00){\makebox(0,0)[cc]{$_a$}}
\put(49.00,75.00){\makebox(0,0)[cc]{$_c$}}
\put(44.00,78.00){\makebox(0,0)[cc]{$_5$}}
\put(36.00,66.00){\makebox(0,0)[cc]{$_7$}}
\put(26.00,66.00){\makebox(0,0)[cc]{$_d$}}
\put(31.00,119.00){\makebox(0,0)[cc]{$_{D_0}$}}
\put(63.00,61.00){\makebox(0,0)[cc]{$_{B_2}$}}
\put(32.00,79.00){\makebox(0,0)[cc]{$_{B_f}$}}
\put(1.00,61.00){\makebox(0,0)[cc]{$_{A_0}$}}
\put(52.00,67.00){\makebox(0,0)[cc]{$_{ab}$}}
\put(82.00,122.00){\makebox(0,0)[cc]{$_{2g}$}}
\put(84.00,127.00){\makebox(0,0)[cc]{$_{46}$}}
\put(75.00,132.00){\makebox(0,0)[cc]{$_{8b}$}}
\put(98.00,170.00){\makebox(0,0)[cc]{$_{1a}$}}
\put(109.00,167.00){\makebox(0,0)[cc]{$_{49}$}}
\put(115.00,170.00){\makebox(0,0)[cc]{$_{27}$}}
\put(127.00,127.00){\makebox(0,0)[cc]{$_{69}$}}
\put(114.00,141.00){\makebox(0,0)[cc]{$_{7g}$}}
\put(109.00,146.00){\makebox(0,0)[cc]{$_{8d}$}}
\put(100.00,142.00){\makebox(0,0)[cc]{$_{a\!f}$}}
\put(131.00,122.00){\makebox(0,0)[cc]{$_{1f}$}}
\put(137.00,132.00){\makebox(0,0)[cc]{$_{bd}$}}
\put(85.00,146.00){\makebox(0,0)[cc]{$_f$}}
\put(91.00,157.00){\makebox(0,0)[cc]{$_d$}}
\put(109.00,161.00){\makebox(0,0)[cc]{$_b$}}
\put(109.00,151.00){\makebox(0,0)[cc]{$_6$}}
\put(122.00,157.00){\makebox(0,0)[cc]{$_8$}}
\put(128.00,146.00){\makebox(0,0)[cc]{$_g$}}
\put(88.00,135.00){\makebox(0,0)[cc]{$_1$}}
\put(94.00,138.00){\makebox(0,0)[cc]{$_9$}}
\put(124.00,135.00){\makebox(0,0)[cc]{$_2$}}
\put(119.00,138.00){\makebox(0,0)[cc]{$_4$}}
\put(111.00,126.00){\makebox(0,0)[cc]{$_7$}}
\put(101.00,126.00){\makebox(0,0)[cc]{$_a$}}
\put(106.00,179.00){\makebox(0,0)[cc]{$_{D_e}$}}
\put(138.00,121.00){\makebox(0,0)[cc]{$_{F_c}$}}
\put(75.00,121.00){\makebox(0,0)[cc]{$_{F_5}$}}
\put(107.00,139.00){\makebox(0,0)[cc]{$_{D_3}$}}
\put(82.00,62.00){\makebox(0,0)[cc]{$_{9g}$}}
\put(84.00,67.00){\makebox(0,0)[cc]{$_{ef}$}}
\put(75.00,72.00){\makebox(0,0)[cc]{$_{3d}$}}
\put(98.00,110.00){\makebox(0,0)[cc]{$_{8c}$}}
\put(109.00,107.00){\makebox(0,0)[cc]{$_{4f}$}}
\put(115.00,110.00){\makebox(0,0)[cc]{$_{5g}$}}
\put(127.00,67.00){\makebox(0,0)[cc]{$_{4e}$}}
\put(113.00,82.00){\makebox(0,0)[cc]{$_{59}$}}
\put(109.00,86.00){\makebox(0,0)[cc]{$_{2d}$}}
\put(99.00,82.00){\makebox(0,0)[cc]{$_{1c}$}}
\put(131.00,62.00){\makebox(0,0)[cc]{$_{18}$}}
\put(137.00,72.00){\makebox(0,0)[cc]{$_{23}$}}
\put(85.00,86.00){\makebox(0,0)[cc]{$_1$}}
\put(91.00,97.00){\makebox(0,0)[cc]{$_2$}}
\put(109.00,101.00){\makebox(0,0)[cc]{$_3$}}
\put(109.00,91.00){\makebox(0,0)[cc]{$_e$}}
\put(122.00,97.00){\makebox(0,0)[cc]{$_d$}}
\put(128.00,86.00){\makebox(0,0)[cc]{$_9$}}
\put(88.00,75.00){\makebox(0,0)[cc]{$_8$}}
\put(94.00,78.00){\makebox(0,0)[cc]{$_4$}}
\put(124.00,75.00){\makebox(0,0)[cc]{$_g$}}
\put(119.00,78.00){\makebox(0,0)[cb]{$_f$}}
\put(111.00,66.00){\makebox(0,0)[cc]{$_5$}}
\put(101.00,66.00){\makebox(0,0)[cc]{$_c$}}
\put(106.00,119.00){\makebox(0,0)[cc]{$_{A_a}$}}
\put(138.00,61.00){\makebox(0,0)[cc]{$_{C_b}$}}
\put(75.00,61.00){\makebox(0,0)[cc]{$_{C_6}$}}
\put(107.00,79.00){\makebox(0,0)[cc]{$_{A_7}$}}
\put(157.00,122.00){\makebox(0,0)[cc]{$_{36}$}}
\put(159.00,127.00){\makebox(0,0)[cc]{$_{be}$}}
\put(150.00,132.00){\makebox(0,0)[cc]{$_{1g}$}}
\put(173.00,170.00){\makebox(0,0)[cc]{$_{2f}$}}
\put(184.00,167.00){\makebox(0,0)[cc]{$_{5b}$}}
\put(190.00,170.00){\makebox(0,0)[cc]{$_{6c}$}}
\put(202.00,127.00){\makebox(0,0)[cc]{$_{5e}$}}
\put(189.00,142.00){\makebox(0,0)[cc]{$_{3c}$}}
\put(186.00,146.00){\makebox(0,0)[cc]{$_{g\infty}$}}
\put(174.00,142.00){\makebox(0,0)[cc]{$_{0f}$}}
\put(206.00,122.00){\makebox(0,0)[cc]{$_{02}$}}
\put(213.00,132.00){\makebox(0,0)[cc]{$_{1\infty}$}}
\put(160.00,146.00){\makebox(0,0)[cc]{$_0$}}
\put(166.00,157.00){\makebox(0,0)[cc]{$_\infty$}}
\put(184.00,161.00){\makebox(0,0)[cc]{$_1$}}
\put(184.00,151.00){\makebox(0,0)[cc]{$_e$}}
\put(197.00,157.00){\makebox(0,0)[cc]{$_g$}}
\put(203.00,146.00){\makebox(0,0)[cc]{$_3$}}
\put(163.00,135.00){\makebox(0,0)[cc]{$_2$}}
\put(169.00,138.00){\makebox(0,0)[cc]{$_5$}}
\put(199.00,135.00){\makebox(0,0)[cc]{$_6$}}
\put(194.00,138.00){\makebox(0,0)[cc]{$_b$}}
\put(186.00,126.00){\makebox(0,0)[cc]{$_c$}}
\put(176.00,126.00){\makebox(0,0)[cc]{$_f$}}
\put(181.00,179.00){\makebox(0,0)[cc]{$_{A_0}$}}
\put(213.00,121.00){\makebox(0,0)[cc]{$_{E_1}$}}
\put(150.00,121.00){\makebox(0,0)[cc]{$_{F_0}$}}
\put(182.00,79.00){\makebox(0,0)[cc]{$_{B_8}$}}
\put(181.00,119.00){\makebox(0,0)[cc]{$_{F_0}$}}
\put(157.00,62.00){\makebox(0,0)[cc]{$_{3a}$}}
\put(159.00,67.00){\makebox(0,0)[cc]{$_{7e}$}}
\put(150.00,72.00){\makebox(0,0)[cc]{$_{89}$}}
\put(173.00,110.00){\makebox(0,0)[cc]{$_{1g}$}}
\put(185.00,107.00){\makebox(0,0)[cc]{$_{be}$}}
\put(191.00,110.00){\makebox(0,0)[cc]{$_{36}$}}
\put(202.00,67.00){\makebox(0,0)[cc]{$_{7b}$}}
\put(189.00,82.00){\makebox(0,0)[cc]{$_{6a}$}}
\put(186.00,86.00){\makebox(0,0)[cc]{$_{8\infty}$}}
\put(174.00,82.00){\makebox(0,0)[cc]{$_{0g}$}}
\put(206.00,62.00){\makebox(0,0)[cc]{$_{01}$}}
\put(213.00,72.00){\makebox(0,0)[cc]{$_{9\infty}$}}
\put(160.00,86.00){\makebox(0,0)[cc]{$_0$}}
\put(166.00,97.00){\makebox(0,0)[cc]{$_\infty$}}
\put(184.00,101.00){\makebox(0,0)[cc]{$_9$}}
\put(184.00,91.00){\makebox(0,0)[cc]{$_7$}}
\put(197.00,97.00){\makebox(0,0)[cc]{$_8$}}
\put(203.00,86.00){\makebox(0,0)[cc]{$_a$}}
\put(163.00,75.00){\makebox(0,0)[cc]{$_1$}}
\put(169.00,78.00){\makebox(0,0)[cc]{$_b$}}
\put(199.00,75.00){\makebox(0,0)[cc]{$_3$}}
\put(194.00,78.00){\makebox(0,0)[cc]{$_e$}}
\put(186.00,66.00){\makebox(0,0)[cc]{$_6$}}
\put(176.00,66.00){\makebox(0,0)[cc]{$_g$}}
\put(213.00,61.00){\makebox(0,0)[cc]{$_{B_9}$}}
\put(150.00,61.00){\makebox(0,0)[cc]{$_{C_0}$}}
\put(182.00,139.00){\makebox(0,0)[cc]{$_{E_g}$}}
\put(76.00,4.00){\circle{2.00}} \put(136.00,4.00){\circle{2.00}}
\put(106.00,54.00){\circle{2.00}} \put(106.00,22.00){\circle{2.00}}
\emline{76.00}{3.00}{41}{136.00}{3.00}{42}
\emline{137.00}{5.00}{43}{107.00}{55.00}{44}
\emline{105.00}{55.00}{45}{75.00}{5.00}{46}
\emline{106.00}{53.00}{47}{106.00}{23.00}{48}
\emline{77.00}{5.00}{49}{105.00}{21.00}{50}
\emline{107.00}{21.00}{51}{135.00}{5.00}{52}
\put(82.00,1.00){\makebox(0,0)[cc]{$_{sy}$}}
\put(84.00,6.00){\makebox(0,0)[cc]{$_{\ell r}$}}
\put(75.00,11.00){\makebox(0,0)[cc]{$_{kx}$}}
\put(98.00,49.00){\makebox(0,0)[cc]{$_{nw}$}}
\put(109.00,46.00){\makebox(0,0)[cc]{$_{rz}$}}
\put(115.00,49.00){\makebox(0,0)[cc]{$_{ms}$}}
\put(128.00,7.00){\makebox(0,0)[cc]{$_{\ell z}$}}
\put(114.00,21.00){\makebox(0,0)[cc]{$_{my}$}}
\put(110.00,25.00){\makebox(0,0)[cc]{$_{px}$}}
\put(99.00,21.00){\makebox(0,0)[cc]{$_{nq}$}}
\put(131.00,1.00){\makebox(0,0)[cc]{$_{qw}$}}
\put(137.00,11.00){\makebox(0,0)[cc]{$_{kp}$}}
\put(85.00,25.00){\makebox(0,0)[cc]{$_q$}}
\put(91.00,36.00){\makebox(0,0)[cc]{$_p$}}
\put(109.00,40.00){\makebox(0,0)[cc]{$_k$}}
\put(109.00,30.00){\makebox(0,0)[cc]{$_\ell$}}
\put(122.00,36.00){\makebox(0,0)[cc]{$_x$}}
\put(128.00,25.00){\makebox(0,0)[cc]{$_y$}}
\put(88.00,14.00){\makebox(0,0)[cc]{$_w$}}
\put(94.00,17.00){\makebox(0,0)[cc]{$_z$}}
\put(124.00,14.00){\makebox(0,0)[cc]{$_s$}}
\put(119.00,17.00){\makebox(0,0)[cc]{$_r$}}
\put(111.00,5.00){\makebox(0,0)[cc]{$_m$}}
\put(101.00,5.00){\makebox(0,0)[cc]{$_n$}}
\put(53.00,167.00){\makebox(0,0)[cc]{$_{A^0}$}}
\put(128.00,167.00){\makebox(0,0)[cc]{$_{B^0}$}}
\put(53.00,107.00){\makebox(0,0)[cc]{$_{D^0}$}}
\put(128.00,107.00){\makebox(0,0)[cc]{$_{E^0}$}}
\put(203.00,107.00){\makebox(0,0)[cc]{$_{F^0}$}}
\put(203.00,167.00){\makebox(0,0)[cc]{$_{C^0}$}}
\put(128.00,46.00){\makebox(0,0)[cc]{$_{\Sigma^0}$}}
\put(181.00,83.00){\circle{2.00}}
\put(104.00,174.00){\vector(1,1){1.00}}
\put(79.00,127.00){\vector(-1,-1){1.00}}
\put(133.00,127.00){\vector(2,-1){2.00}}
\put(104.00,114.00){\vector(1,1){1.00}}
\put(135.00,69.00){\vector(1,-2){1.00}}
\emline{31.00}{144.00}{53}{31.00}{174.00}{54}
\emline{106.00}{144.00}{55}{106.00}{174.00}{56}
\emline{31.00}{84.00}{57}{31.00}{114.00}{58}
\emline{106.00}{84.00}{59}{106.00}{114.00}{60}
\emline{76.00}{124.00}{61}{136.00}{124.00}{62}
\emline{76.00}{64.00}{63}{136.00}{64.00}{64}
\emline{30.00}{116.00}{65}{0.00}{66.00}{66}
\emline{2.00}{66.00}{67}{30.00}{82.00}{68}
\emline{62.00}{126.00}{69}{32.00}{176.00}{70}
\emline{32.00}{142.00}{71}{60.00}{126.00}{72}
\put(179.00,174.00){\vector(1,1){1.00}}
\put(184.00,113.00){\vector(-1,1){2.00}}
\emline{212.00}{126.00}{73}{182.00}{176.00}{74}
\emline{182.00}{142.00}{75}{210.00}{126.00}{76}
\emline{180.00}{116.00}{77}{150.00}{66.00}{78}
\emline{152.00}{66.00}{79}{180.00}{82.00}{80}
\emline{181.00}{144.00}{81}{181.00}{174.00}{82}
\emline{181.00}{84.00}{83}{181.00}{114.00}{84}
\put(29.00,174.00){\vector(1,1){1.00}}
\put(109.00,173.00){\vector(-1,1){2.00}}
\put(34.00,173.00){\vector(-1,1){2.00}}
\put(29.00,114.00){\vector(1,1){1.00}}
\put(34.00,113.00){\vector(-1,1){2.00}}
\put(109.00,113.00){\vector(-1,1){2.00}}
\put(79.00,67.00){\vector(-1,-1){1.00}}
\put(133.00,67.00){\vector(2,-1){2.00}}
\put(31.00,146.00){\vector(0,-1){2.00}}
\put(181.00,146.00){\vector(0,-1){2.00}}
\put(31.00,86.00){\vector(0,-1){2.00}}
\put(181.00,86.00){\vector(0,-1){2.00}}
\put(184.00,173.00){\vector(-1,1){2.00}}
\put(179.00,114.00){\vector(1,1){1.00}}
\end{picture}
%\vspace{-2mm}
\caption{Symmetry of edge labels in copies of $K_4$ in
$\mathcal Y$, for $i=0$}
\end{figure}
\noindent Figure 3 illustrates the left side of (5) for $i=0$ in
terms of edge labels, where edges of $\mathcal Y$ arising from pairs
of 3-arcs of $\mathcal S$ identically (resp. oppositely) fastened
according to (1) are shown oriented (resp. unoriented) accordingly.
Observe the edges oriented in
$$\begin{array}{lll}
^{A^0\,:\,D_0C_0,\;C_0E_4,\;C_0E_d;}
_{D^0\,:\,A_0D_0,\;B_2D_0,\;D_0B_f;}&
^{B^0\,:\,D_3F_5,\;F_5D_e,\;D_3F_c,\,F_cD_e;}
_{E^0\,:\,A_7C_6,\;A_7C_b,\;C_bA_a,\,C_6A_a;}&
^{C^0\,:\,F_0A_0,\;A_0E_1,\;E_gA_0;}
_{F^0\,:\,F_0B_8,\;B_9F_0,\;C_0F_0.}\end{array}$$ By uniformly
adding successively $1\in\Z_{17}$, each of these 6 cases yields 16
additional ones. This yields the 102 edge-labeled copies of $K_4$ in
$\mathcal Y$. If the two points of $PG(1,17)$ labeling near its
center each edge $\epsilon$ in the figure are disposed as shown,
labeling the respective flags of $\epsilon$, then the 6 cases may be
indicated uniquely as $(kl,mn)(pq,rs)(xy,zw)$, where the position of
the labels $k,\ell,m,n,p,q,r,s,w,x,y,z$ is as in the referential
depiction $\Sigma^0$ of a copy of $K_4$ in the lower part of the
figure. Then, the flag-label triples at the upper, middle,
lower-right and lower-left vertices of this depiction are
respectively $kpx$, $\ell rz$, $msy$ and $nqw$. Moreover, the 6
points of $PG(1,17)$ in each of these copies of $K_4$ not
participating of its edge labeling conform a unique sextet $\chi$
which is not a vertex of $\mathcal S$ as characterized in Section 2.
However, $\chi$ is a sextet of an alternative labeling of $\mathcal
S$ happening via the remaining 102 sextets (of the total of 204).
These 102 alternative sextets are the images of the 102 vertices of
$\mathcal S$ via multiplication of indices in $PG(1,17)$ times $3\in
GF(17)$, operation that coincides with the duality $\phi$ expressed
in (6) above. This proves the assertion in Theorem 4.1 below that
the vertices and copies of $K_4$ of $\mathcal S$ are the points and
lines of a self-dual 1-configuration $(102_4)_1$, which in turn has
$\mathcal Y$ as its Menger graph. Correspondingly, the vertex labels
in $\Sigma^i$ are the sextets $(rz,ms,nw)$, $(px,nq$, $my)$
$(kp,\ell z,qw)$ and $(kx,\ell r,sy)$.

\noindent A procedure that allows to determine which point of
$PG(1,17)$ labels which flag in a copy of $K_4$ as in Figure 3 is
given as follows:

{\bf(i)} A triangle $\Delta$ in a copy $\nabla$ of
$K_4$ in $\mathcal Y$, say $\Delta=(C_0E_4D_0)$ in $\nabla=A^0$,
arises from a 9-cycle $\Pi^j=(\Pi_0^j\ldots\Pi_8^j)$ in $\mathcal S$
with associated permutation
$\pi^j=(\pi_0^j\ldots\pi_8^j)(\xi_0^j\ldots\xi_8^j)$ as displayed in
Section 2, in this case $\Pi^j=Y^2$ with $\pi^j=x^2$; and

{\bf(ii)} by labeling each edge $\Pi_i^j\Pi_{i+1}^j$ of $\Pi^j$ just
by $\pi_i^j$, it holds that the flag label of edge
$\epsilon=\Pi_i^j\Pi_{i+3}^j$ at $\Pi_i^j$ is $\pi_{i+1}^j$ while
the flag label of $\epsilon$ at $\Pi_{i+3}^j$ is $\pi_{i+5}^j$,
where $i=0,3,6$.

\noindent The distance-3 digraphs of the directed 9-cycles $\Pi^0$
of $\mathcal S$ are composed by the following triples of disjoint
directed triangles of $\mathcal Y$:
$$\begin{array}{l}
^{S^0\,\rightarrow\,\{D^0\setminus D_0=(B_2A_0B_f),\,\,\,E^9\setminus C_3=(A_2A_gC_f),\,\,\,E^8\setminus C_8=(A_1A_fC_2)\};}
_{T^0\,\rightarrow\,\{A^0\setminus C_0=(E_dD_0E_4),\,\,\,B^g\setminus F_b=(D_dD_2F_4),\,\,\,B^1\setminus F_1=(D_fD_4F_d)
\};}\\
^{U^0\,\rightarrow\,\{F^0\setminus F_0=(B_9C_0B_8),\,\,\,E^f\setminus A_5=(C_9C_4A_8),\,\,\,E^2\setminus A_2=(C_dC_8A_9)\};}
_{V^0\,\rightarrow\,\{C^0\setminus A_0=(E_gF_0E_1),\,\,\,B^4\setminus D_7=(F_gF_9D_1),\,\,\,B^d\setminus D_d=(F_8F_1D_g)\};}\\
^{W^0\rightarrow\{F^0\setminus C_0=(B_9F_0B_8),\,\,\,C^8\setminus E_7=(E_9F_8A_8),\,\,\,\,C^9\setminus C_9=(F_9E_8A_9)\};}
_{X^0\,\rightarrow\{C^0\setminus F_0=(E_gA_0E_1),\,\,\,D^1\setminus B_3=(B_gA_1D_1),\,\,D^g\setminus D_g=(A_gB_1D_g)\};}\\
^{Y^0\,\rightarrow\{D^0\setminus A_0=(B_2D_0B_f),\,\,A^f\setminus E_b=(E_2D_fC_f),\,\,A^2\setminus A_2=(D_2E_fC_2)\};}
_{Z^0\,\rightarrow\{A^0\setminus D_0=(E_dC_0E_4),\,\,\,F^4\setminus D_c=(B_dC_4F_4),\,\,\,\,F^d\setminus F_d=(C_dB_4F_d)\}.}
\end{array}$$
\noindent This way, it can be seen that $\mathcal Y$ is a $K_4$-UH
graph. However, in view of Beineke's characterization of line graphs
\cite{Bei} and observing that $\mathcal Y$ contains induced copies
of $K_{1,3}$\,, which are forbidden for line graphs of simple
graphs, we conclude that $\mathcal Y$ is non-line-graphical.

\begin{thm} $\mathcal Y$ is both the Menger graph of a $K_4$-UH self-dual
$1$-con\-fi\-gu\-ra\-tion $(102_4)_1$ and a non-line-graphical
$\{K_4\}_{102}^{4}$-graph. Moreover, $\mathcal Y$ is
arc-tran\-si\-tive with regular degree $12$, diameter $3$, distance
distribution $(1$, $12$, $78$, $11)$ and automorphism group
$PSL(2,17)$ of order $2448$. Its associated Levi graph is a
$2$-arc-transitive graph with regular degree $4$, diameter $6$,
distance distribution $(1,4,12,36,78,62,11)$ and automorphism group
$SL(2,17)$ of order $4896$.
\end{thm}

\proof It remains to prove that $\mathcal Y$ is $K_4$-UH, which uses
(1) and more specifically (4) above. In fact, consider an
isomorphism $\Psi:\Theta_1\rightarrow\Theta_2$ between copies
$\Theta_1,\Theta_2$ of $K_4$ in $\mathcal Y$. Each $\Theta_i$,
($i=1,2$), arises from 4 9-cycles $\gamma_9=\theta_i^j$ in $\mathcal
S$, ($j=1,2,3,4$), whose union is a subgraph $\overline{\Theta}_i$
of $\mathcal S$ with 4 vertices $v_i^j$ of degree 3 and 12 vertices
of degree 2 that are the internal vertices of 6 3-paths $P_4$ whose
ends are the vertices $v_i^j$. For example, the vertices
$v_1^1=B_0,v_1^2=B_1,v_1^3=F_9,v_1^4=C_9,v_2^1=B_1,v_2^2=B_2,v_2^3=F_a,v_2^4=C_a$
in $\mathcal S$ determine such subgraphs $\Theta_1,\Theta_2$ in
$\mathcal Y$ and $\overline{\Theta}_1,\overline{\Theta}_2$ in
$\mathcal S$. Clearly, $\Psi$ induces an isomorphism
$\overline{\Psi}:\overline{\Theta}_1\rightarrow\overline{\Theta}_2$
that sends say each $v_1^j$ onto its corresponding $v_2^j$,
($j=1,2,3,4$). As an automorphism $\overline{\overline{\Psi}}$ of
$\mathcal S$ exists that extends $\overline{\Psi}$, then
$\overline{\overline{\Psi}}$ determines an automorphism of $\mathcal
Y$ that restricts to $\Psi$, showing that $\mathcal Y$ is a $K_4$-UH
graph.\qfd

\section{Definitions to deal with the copies of $L(Q_3)$}

\noindent If $H$ is a graph with an edge partition
$\Omega=\Omega(H)$ into 2-paths, then a graph $G$ is $\Omega$-{\it
preserving} $H$-UH if every $\Omega$-preserving isomorphism between
two induced copies of $H$ in $G$ extends to an automorphism of $G$.
If $M$ is a subgraph of $H$ and if $G$ is both $M$-UH, and
$\Omega$-preserving $H$-UH, then $G$ is an $\Omega$-{\it preserving}
$\{H\}_{M}$-{\it UH graph} if, for each induced copy $H_0$ of $H$ in
$G$ containing an induced copy $M_0$ of $M$, there is just one
induced copy $H_1\neq H_0$ of $H$ in $G$ such that:

{\bf(a)} $V(H_0)\cap V(H_1)=V(M_0)$;

{\bf(b)} $E(H_0)\cap E(H_1)=E(M_0)$; and

{\bf(c)} the edges of $M_0$ are in distinct 2-paths both in
$\Omega(H_0)$ and $\Omega(H_1)$.

\noindent A graph $G$ is $rK_s$-{\it frequent} if every edge $e$ of
$G$ is intersection of exactly $r$ induced copies of $K_s$, these
copies having only $e$ and its ends in common. For example, $K_4$ is
$2K_3$-frequent and $L(Q_3)$ is $1K_3$-frequent. A graph $G$ is
$\{H_2,H_1\}_{K_3}$-{\it UH}, where $H_i$ is $iK_3$-frequent
($i=1,2$) if:

{\bf(d)} $G$ is $H_2$-UH and edge-disjoint union of induced copies
of $H_2$;

{\bf(e)} there is a partition $\Omega$ of $H_1$ into 2-paths and $G$
is $\Omega$-preserving \hspace*{1cm} $\{H_1\}_{K_3}$-UH; and

{\bf(f)} each induced copy of $H_2$ in $G$ has each induced copy of
$K_3$ in \hspace*{1cm} common with exactly two induced copies of
$H_1$ in $G$.

\noindent Theorem 6.1 shows that $\mathcal Y$ is
$\{K_4,L(Q_3)\}_{K_3}$-UH. This allows to gather information on
${\mathcal S}_2$ and ${\mathcal S}_4$, leading to ${\mathcal
Y}={\mathcal S}_3$ in Theorem 7.1.

\section{The $K_4$-UH graph $\mathcal Y$ is $\{K_4,L(Q_3)\}_{K_3}$-UH}

\begin{figure}[htp]
\vspace*{-6mm} \hspace*{.5mm}
  \includegraphics[scale=0.33]{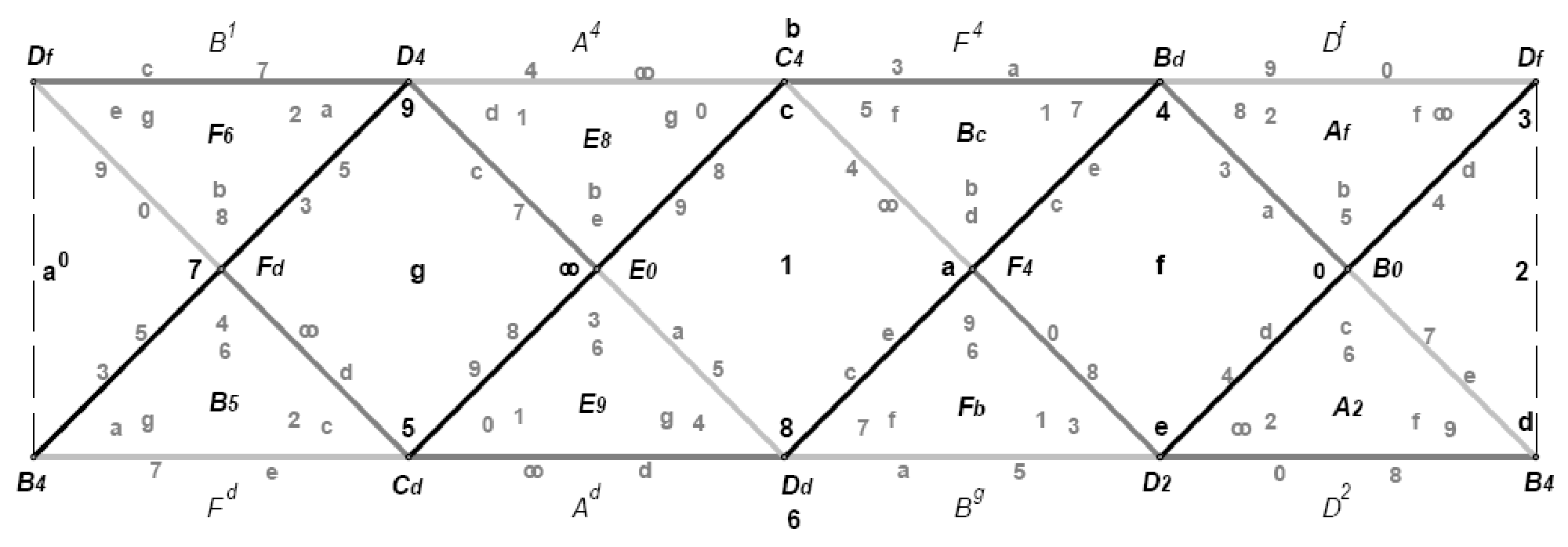}%0.41]{itacolor.eps}\\
\caption{Toroidal cutout representation of $a^0$}
\end{figure}

\noindent Recall from Section 4 that each copy of $K_4$ in $\mathcal
Y$ arises from the distance-3 digraphs of 4 directed 9-cycles of
$\mathcal S$. The subgraph of $\mathcal S$ spanned by these 4
9-cycles contains 4 degree-3 vertices (which are tails and heads of
corresponding 3-arcs) and 12 degree-2 vertices (internal vertices of
those 3-arcs). These 12 vertices induce a copy $\mathcal L$ of
$L(Q_3)$ in $\mathcal Y$. For the copy $A^0$ of $K_4$ in $\mathcal
Y$, the corresponding copy ${\mathcal L}=a^0$ of $L(Q_3)$ in
$\mathcal Y$ can be represented as in the big rectangle $\mathcal R$
in Figure 4, where:

{\bf(a)} the leftmost and rightmost dashed lines of $\mathcal R$ are
to be identified by parallel translation;

{\bf(b)} each of the 8 shown triangles $\Delta$ forms part of a
corresponding copy $\mathcal r$ of $K_4$ cited on the exterior of
$\mathcal R$ about the horizontal edge of $\Delta$, while its 4th
vertex is cited at the center of $\Delta$; and

{\bf(c)} the edges are colored via a partition $\Omega$ into 2-paths
$P_3$, the edges of each $P_3$ with a common color from a set of 3
colors: {\bf(i)} black; {\bf(ii)} light-gray; {\bf(iii)} dark-gray;
the 3 colors are present together in every triangle, and opposite
edges in every induced 4-cycle, or 4-{\it hole}, have a common
color, a total of two colors per 4-hole.

\noindent For $\sigma=a,b,c,d,e,f$, the copies $\sigma^0$ of
$L(Q_3)$ are expressed by means of the data contained in Figure 4 as
follows:

$$\begin{array}{l}
^{a^0:(D_fD_4C_4B_d)\!(B_4C_dD_dD_2)F_dE_0F_4B_0\!(B^1F_6A^4E_8F^4B_cD^fA_f)\!(F^dB_5A^dE_9B^gF_bD^2A_2)}
_{b^0:(D_5D_gE_cF_d)\!(D_cF_4E_5D_1)F_eE_3E_eF_3\!(B^2F_7A^gC_gC^dA_dB^8D_b)\!(B^9D_6C^4A_4A^1C_1B^aD_f)}\\
^{c^0:(F_8F_1A_1B_g)\!(B_1A_gF_gF_9)D_gE_0D_1B_0\!(B^dD_aC^1E_2D^1B_3F^8C_8)\!(D^gB_eC^gE_fB^4D_7F^9C_9)}
_{d^0:(A_1A_fD_fE_2)\!(E_fD_2A_2A_g)C_2B_0C_fE_0\!(E^8C_eD^fB_dA^fE_bC^1F_1)\!(A^2E_6D^2B_4E^9C_3C^gF_g)}\\
^{e^0:(A_6A_9B_bC_2)\!(A_bC_fB_6A_8)C_aB_7B_aC_7\!(E^gC_5D^9D_9F^2F_2E^dA_3)\!(E^4A_eF^fF_fD^aE_8E^1C_c)}
_{f^0:(C_4C_9F_9E_8)\!(E_9F_8C_8C_d)A_8B_0A_9E_0\!(E^fA_5F^9B_1C^9E_aA^4D_4)\!(C^8E_7F^8B_gE^2A_cA^dD_d)}
\end{array}$$
and  their translations mod 17 are denoted $\sigma^i$,  for $0\ne
i\in\Z_{17}$ (uniformly translating all involved subscripts and
superscripts). Each copy $\sigma^i$ of $L(Q_3)$ admits an edge
partition $\Omega=\Omega(\sigma^i)$ into $j$-colored 2-paths
($j\in\{1,2,3\}$) so that each (monochromatic) 2-path in an
$\Omega(\sigma^i)$ is shared only by one other copy of $L(Q_3)$ in
$\mathcal Y$ (as in Theorem 6.1{\bf(3)}, below). We may write
\begin{equation}\sigma^i=\sigma^i_1\cup \sigma^i_2\cup
\sigma^i_3,\end{equation} \noindent to stress the color partition of
$\sigma^i$ into its black, light-gray and dark-gray subgraphs, which
are copies of the disconnected graph $4P_3$ (formed by 4 disjoint
copies of $P_3$) as in Figure 4 for $\sigma^i=a^0$. The edge labels
of $\sigma^0$ in Figure 4 (shown in gray type) and of all the other
$\sigma^i$\thinspace s are taken as the flag labels for
$i=0,\ldots,g$ in Figure 3. The relation and location of these flag
labels justifies a labeling of the 12 vertices and 6 4-holes as
shown with symbols $0,\ldots,g,\infty$ (in black type) in Figure 4,
the sole edge-label notation to be used ahead.

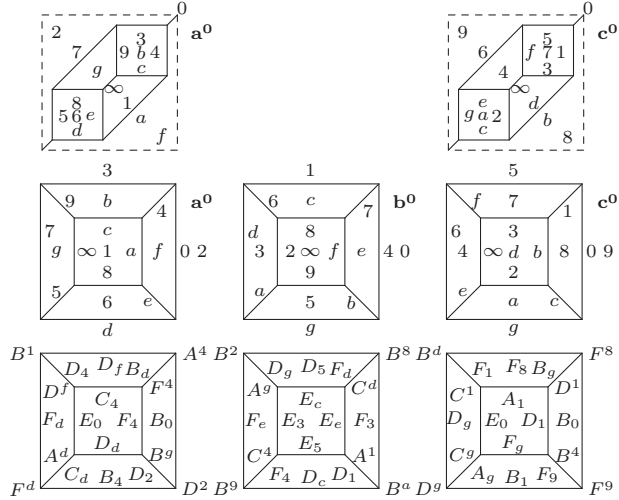
\begin{figure}[htp]
\hspace*{1.2cm}
\unitlength=0.45mm \special{em:linewidth 0.4pt}
\linethickness{0.4pt}
\begin{picture}(181.50,143.00)
\put(33.50,27.00){\makebox(0,0)[cc]{$_{C_4}$}}
\put(33.50,14.00){\makebox(0,0)[cc]{$_{D_d}$}}
\put(39.50,21.00){\makebox(0,0)[cc]{$_{F_4}$}}
\put(28.50,21.00){\makebox(0,0)[cc]{$_{E_0}$}}
\put(42.50,36.00){\makebox(0,0)[cc]{$_{B_d}$}}
\put(24.50,5.00){\makebox(0,0)[cc]{$_{C_d}$}}
\put(24.50,36.00){\makebox(0,0)[cc]{$_{D_4}$}}
\put(49.50,21.00){\makebox(0,0)[cc]{$_{B_0}$}}
\put(43.50,5.00){\makebox(0,0)[cc]{$_{D_2}$}}
\put(34.50,4.00){\makebox(0,0)[cc]{$_{B_4}$}}
\put(34.50,37.00){\makebox(0,0)[cc]{$_{D_f}$}}
\put(17.50,21.00){\makebox(0,0)[cc]{$_{F_d}$}}
\put(33.50,71.00){\makebox(0,0)[cc]{$_1$}}
\put(33.50,77.00){\makebox(0,0)[cc]{$_c$}}
\put(33.50,65.00){\makebox(0,0)[cc]{$_8$}}
\put(40.50,71.00){\makebox(0,0)[cc]{$_a$}}
\put(27.50,71.00){\makebox(0,0)[cc]{$_\infty$}}
\put(48.50,71.00){\makebox(0,0)[cc]{$_f$}}
\put(49.50,83.00){\makebox(0,0)[cc]{$_4$}}
\put(33.50,56.00){\makebox(0,0)[cc]{$_6$}}
\put(18.50,59.00){\makebox(0,0)[cc]{$_5$}}
\put(18.50,71.00){\makebox(0,0)[cc]{$_g$}}
\put(22.50,86.00){\makebox(0,0)[cc]{$_9$}}
\put(33.50,86.00){\makebox(0,0)[cc]{$_b$}}
\put(56.50,71.00){\makebox(0,0)[cc]{$_0$}}
\put(61.50,71.00){\makebox(0,0)[cc]{$_2$}}
\put(61.50,86.00){\makebox(0,0)[cc]{$\mathbf{_{a^0}}$}}
\put(45.50,56.00){\makebox(0,0)[cc]{$_e$}}
\put(93.50,71.00){\makebox(0,0)[cc]{$_\infty$}}
\put(93.50,77.00){\makebox(0,0)[cc]{$_8$}}
\put(93.50,65.00){\makebox(0,0)[cc]{$_9$}}
\put(100.50,71.00){\makebox(0,0)[cc]{$_f$}}
\put(87.50,71.00){\makebox(0,0)[cc]{$_2$}}
\put(108.50,71.00){\makebox(0,0)[cc]{$_e$}}
\put(110.50,83.00){\makebox(0,0)[cc]{$_7$}}
\put(93.50,56.00){\makebox(0,0)[cc]{$_5$}}
\put(78.50,59.00){\makebox(0,0)[cc]{$_a$}}
\put(78.50,71.00){\makebox(0,0)[cc]{$_3$}}
\put(82.50,86.00){\makebox(0,0)[cc]{$_6$}}
\put(93.50,86.00){\makebox(0,0)[cc]{$_c$}}
\put(116.50,71.00){\makebox(0,0)[cc]{$_4$}}
\put(121.50,71.00){\makebox(0,0)[cc]{$_0$}}
\put(121.50,86.00){\makebox(0,0)[cc]{$_\mathbf{{b^0}}$}}
\put(105.50,56.00){\makebox(0,0)[cc]{$_b$}}
\put(33.50,48.00){\makebox(0,0)[cc]{$_d$}}
\put(93.50,48.00){\makebox(0,0)[cc]{$_g$}}
\put(33.50,95.00){\makebox(0,0)[cc]{$_3$}}
\put(93.50,95.00){\makebox(0,0)[cc]{$_1$}}
\put(16.50,77.00){\makebox(0,0)[cc]{$_7$}}
\put(76.50,77.00){\makebox(0,0)[cc]{$_d$}}
\put(153.50,71.00){\makebox(0,0)[cc]{$_d$}}
\put(153.50,77.00){\makebox(0,0)[cc]{$_3$}}
\put(153.50,65.00){\makebox(0,0)[cc]{$_2$}}
\put(160.50,71.00){\makebox(0,0)[cc]{$_b$}}
\put(147.50,71.00){\makebox(0,0)[cc]{$_\infty$}}
\put(168.50,71.00){\makebox(0,0)[cc]{$_8$}}
\put(169.50,83.00){\makebox(0,0)[cc]{$_1$}}
\put(153.50,56.00){\makebox(0,0)[cc]{$_a$}}
\put(138.50,59.00){\makebox(0,0)[cc]{$_e$}}
\put(138.50,71.00){\makebox(0,0)[cc]{$_4$}}
\put(142.50,86.00){\makebox(0,0)[cc]{$_f$}}
\put(153.50,86.00){\makebox(0,0)[cc]{$_7$}}
\put(176.50,71.00){\makebox(0,0)[cc]{$_0$}}
\put(181.50,71.00){\makebox(0,0)[cc]{$_9$}}
\put(181.50,86.00){\makebox(0,0)[cc]{$\mathbf{_{c^0}}$}}
\put(165.50,56.00){\makebox(0,0)[cc]{$_c$}}
\put(153.50,48.00){\makebox(0,0)[cc]{$_g$}}
\put(153.50,95.00){\makebox(0,0)[cc]{$_5$}}
\put(136.50,77.00){\makebox(0,0)[cc]{$_6$}}
\put(93.50,27.00){\makebox(0,0)[cc]{$_{E_c}$}}
\put(93.50,14.00){\makebox(0,0)[cc]{$_{E_5}$}}
\put(99.50,21.00){\makebox(0,0)[cc]{$_{E_e}$}}
\put(88.50,21.00){\makebox(0,0)[cc]{$_{E_3}$}}
\put(102.50,36.00){\makebox(0,0)[cc]{$_{F_d}$}}
\put(84.50,5.00){\makebox(0,0)[cc]{$_{F_4}$}}
\put(84.50,36.00){\makebox(0,0)[cc]{$_{D_g}$}}
\put(109.50,21.00){\makebox(0,0)[cc]{$_{F_3}$}}
\put(103.50,5.00){\makebox(0,0)[cc]{$_{D_1}$}}
\put(94.50,4.00){\makebox(0,0)[cc]{$_{D_c}$}}
\put(94.50,37.00){\makebox(0,0)[cc]{$_{D_5}$}}
\put(77.50,21.00){\makebox(0,0)[cc]{$_{F_e}$}}
\put(153.50,27.00){\makebox(0,0)[cc]{$_{A_1}$}}
\put(153.50,14.00){\makebox(0,0)[cc]{$_{F_g}$}}
\put(159.50,21.00){\makebox(0,0)[cc]{$_{D_1}$}}
\put(148.50,21.00){\makebox(0,0)[cc]{$_{E_0}$}}
\put(162.50,36.00){\makebox(0,0)[cc]{$_{B_g}$}}
\put(144.50,5.00){\makebox(0,0)[cc]{$_{A_g}$}}
\put(144.50,36.00){\makebox(0,0)[cc]{$_{F_1}$}}
\put(169.50,21.00){\makebox(0,0)[cc]{$_{B_0}$}}
\put(163.50,5.00){\makebox(0,0)[cc]{$_{F_9}$}}
\put(154.50,4.00){\makebox(0,0)[cc]{$_{B_1}$}}
\put(154.50,37.00){\makebox(0,0)[cc]{$_{F_8}$}}
\put(137.50,21.00){\makebox(0,0)[cc]{$_{D_g}$}}
\put(179.50,41.00){\makebox(0,0)[cc]{$_{F^8}$}}
\put(179.50,1.00){\makebox(0,0)[cc]{$_{F^9}$}}
\put(128.50,41.00){\makebox(0,0)[cc]{$_{B^d}$}}
\put(128.50,1.00){\makebox(0,0)[cc]{$_{D^g}$}}
\put(138.50,29.00){\makebox(0,0)[cc]{$_{C^1}$}}
\put(138.50,11.00){\makebox(0,0)[cc]{$_{C^g}$}}
\put(169.50,31.00){\makebox(0,0)[cc]{$_{D^1}$}}
\put(169.50,11.00){\makebox(0,0)[cc]{$_{B^4}$}}
\put(119.50,41.00){\makebox(0,0)[cc]{$_{B^8}$}}
\put(119.50,1.00){\makebox(0,0)[cc]{$_{B^a}$}}
\put(68.50,41.00){\makebox(0,0)[cc]{$_{B^2}$}}
\put(68.50,1.00){\makebox(0,0)[cc]{$_{B^9}$}}
\put(78.50,30.00){\makebox(0,0)[cc]{$_{A^g}$}}
\put(78.50,11.00){\makebox(0,0)[cc]{$_{C^4}$}}
\put(109.50,31.00){\makebox(0,0)[cc]{$_{C^d}$}}
\put(109.50,11.00){\makebox(0,0)[cc]{$_{A^1}$}}
\put(59.50,41.00){\makebox(0,0)[cc]{$_{A^4}$}}
\put(59.50,1.00){\makebox(0,0)[cc]{$_{D^2}$}}
\put(8.50,41.00){\makebox(0,0)[cc]{$_{B^1}$}}
\put(8.50,1.00){\makebox(0,0)[cc]{$_{F^d}$}}
\put(18.50,30.00){\makebox(0,0)[cc]{$_{D^f}$}}
\put(18.50,11.00){\makebox(0,0)[cc]{$_{A^d}$}}
\put(49.50,31.00){\makebox(0,0)[cc]{$_{F^4}$}}
\put(49.50,11.00){\makebox(0,0)[cc]{$_{B^g}$}}
\emline{13.50}{51.00}{1}{13.50}{91.00}{2}
\emline{13.50}{91.00}{3}{53.50}{91.00}{4}
\emline{53.50}{91.00}{5}{43.50}{81.00}{6}
\emline{43.50}{81.00}{7}{23.50}{81.00}{8}
\emline{43.50}{81.00}{9}{43.50}{61.00}{10}
\emline{13.50}{91.00}{11}{23.50}{81.00}{12}
\emline{13.50}{51.00}{13}{53.50}{51.00}{14}
\emline{53.50}{51.00}{15}{53.50}{91.00}{16}
\emline{23.50}{81.00}{17}{23.50}{61.00}{18}
\emline{43.50}{61.00}{19}{23.50}{61.00}{20}
\emline{43.50}{61.00}{21}{53.50}{51.00}{22}
\emline{23.50}{61.00}{23}{13.50}{51.00}{24}
\emline{13.50}{1.00}{25}{13.50}{41.00}{26}
\emline{13.50}{41.00}{27}{53.50}{41.00}{28}
\emline{53.50}{41.00}{29}{43.50}{31.00}{30}
\emline{43.50}{31.00}{31}{23.50}{31.00}{32}
\emline{43.50}{31.00}{33}{43.50}{11.00}{34}
\emline{13.50}{41.00}{35}{23.50}{31.00}{36}
\emline{13.50}{1.00}{37}{53.50}{1.00}{38}
\emline{53.50}{1.00}{39}{53.50}{41.00}{40}
\emline{23.50}{31.00}{41}{23.50}{11.00}{42}
\emline{43.50}{11.00}{43}{23.50}{11.00}{44}
\emline{43.50}{11.00}{45}{53.50}{1.00}{46}
\emline{23.50}{11.00}{47}{13.50}{1.00}{48}
\emline{73.50}{91.00}{49}{113.50}{91.00}{50}
\emline{113.50}{91.00}{51}{103.50}{81.00}{52}
\emline{103.50}{81.00}{53}{83.50}{81.00}{54}
\emline{103.50}{81.00}{55}{103.50}{61.00}{56}
\emline{73.50}{91.00}{57}{83.50}{81.00}{58}
\emline{73.50}{51.00}{59}{113.50}{51.00}{60}
\emline{113.50}{51.00}{61}{113.50}{91.00}{62}
\emline{83.50}{81.00}{63}{83.50}{61.00}{64}
\emline{103.50}{61.00}{65}{83.50}{61.00}{66}
\emline{103.50}{61.00}{67}{113.50}{51.00}{68}
\emline{83.50}{61.00}{69}{73.50}{51.00}{70}
\emline{73.50}{1.00}{71}{73.50}{41.00}{72}
\emline{73.50}{41.00}{73}{113.50}{41.00}{74}
\emline{113.50}{41.00}{75}{103.50}{31.00}{76}
\emline{103.50}{31.00}{77}{83.50}{31.00}{78}
\emline{103.50}{31.00}{79}{103.50}{11.00}{80}
\emline{73.50}{41.00}{81}{83.50}{31.00}{82}
\emline{73.50}{1.00}{83}{113.50}{1.00}{84}
\emline{113.50}{1.00}{85}{113.50}{41.00}{86}
\emline{83.50}{31.00}{87}{83.50}{11.00}{88}
\emline{103.50}{11.00}{89}{83.50}{11.00}{90}
\emline{103.50}{11.00}{91}{113.50}{1.00}{92}
\emline{83.50}{11.00}{93}{73.50}{1.00}{94}
\emline{73.50}{51.00}{95}{73.50}{91.00}{96}
\emline{133.50}{51.00}{97}{133.50}{91.00}{98}
\emline{133.50}{91.00}{99}{173.50}{91.00}{100}
\emline{173.50}{91.00}{101}{163.50}{81.00}{102}
\emline{163.50}{81.00}{103}{143.50}{81.00}{104}
\emline{163.50}{81.00}{105}{163.50}{61.00}{106}
\emline{133.50}{91.00}{107}{143.50}{81.00}{108}
\emline{133.50}{51.00}{109}{173.50}{51.00}{110}
\emline{173.50}{51.00}{111}{173.50}{91.00}{112}
\emline{143.50}{81.00}{113}{143.50}{61.00}{114}
\emline{163.50}{61.00}{115}{143.50}{61.00}{116}
\emline{163.50}{61.00}{117}{173.50}{51.00}{118}
\emline{143.50}{61.00}{119}{133.50}{51.00}{120}
\emline{133.50}{1.00}{121}{133.50}{41.00}{122}
\emline{133.50}{41.00}{123}{173.50}{41.00}{124}
\emline{173.50}{41.00}{125}{163.50}{31.00}{126}
\emline{163.50}{31.00}{127}{143.50}{31.00}{128}
\emline{163.50}{31.00}{129}{163.50}{11.00}{130}
\emline{133.50}{41.00}{131}{143.50}{31.00}{132}
\emline{133.50}{1.00}{133}{173.50}{1.00}{134}
\emline{173.50}{1.00}{135}{173.50}{41.00}{136}
\emline{143.50}{31.00}{137}{143.50}{11.00}{138}
\emline{163.50}{11.00}{139}{143.50}{11.00}{140}
\emline{163.50}{11.00}{141}{173.50}{1.00}{142}
\emline{143.50}{11.00}{143}{133.50}{1.00}{144}
\emline{17.00}{104.00}{145}{32.00}{104.00}{146}
\emline{32.00}{104.00}{147}{32.00}{119.00}{148}
\emline{32.00}{119.00}{149}{17.00}{119.00}{150}
\emline{17.00}{119.00}{151}{17.00}{104.00}{152}
\emline{36.00}{123.00}{153}{51.00}{123.00}{154}
\emline{51.00}{123.00}{155}{51.00}{138.00}{156}
\emline{51.00}{138.00}{157}{36.00}{138.00}{158}
\emline{36.00}{138.00}{159}{36.00}{123.00}{160}
\emline{17.00}{119.00}{161}{36.00}{138.00}{162}
\emline{32.00}{119.00}{163}{36.00}{123.00}{164}
\emline{32.00}{104.00}{165}{51.00}{123.00}{166}
\emline{51.00}{138.00}{167}{54.00}{141.00}{168}
\emline{17.00}{104.00}{169}{14.00}{101.00}{170}
\emline{16.00}{101.00}{171}{18.00}{101.00}{172}
\emline{20.00}{101.00}{173}{22.00}{101.00}{174}
\emline{24.00}{101.00}{175}{26.00}{101.00}{176}
\emline{28.00}{101.00}{177}{30.00}{101.00}{178}
\emline{32.00}{101.00}{179}{34.00}{101.00}{180}
\emline{36.00}{101.00}{181}{38.00}{101.00}{182}
\emline{40.00}{101.00}{183}{42.00}{101.00}{184}
\emline{44.00}{101.00}{185}{46.00}{101.00}{186}
\emline{48.00}{101.00}{187}{50.00}{101.00}{188}
\emline{52.00}{101.00}{189}{54.00}{101.00}{190}
\emline{54.00}{103.00}{191}{54.00}{105.00}{192}
\emline{54.00}{107.00}{193}{54.00}{109.00}{194}
\emline{54.00}{111.00}{195}{54.00}{113.00}{196}
\emline{54.00}{115.00}{197}{54.00}{117.00}{198}
\emline{54.00}{119.00}{199}{54.00}{121.00}{200}
\emline{54.00}{123.00}{201}{54.00}{125.00}{202}
\emline{54.00}{127.00}{203}{54.00}{129.00}{204}
\emline{54.00}{131.00}{205}{54.00}{133.00}{206}
\emline{54.00}{135.00}{207}{54.00}{137.00}{208}
\emline{54.00}{139.00}{209}{54.00}{141.00}{210}
\emline{52.00}{141.00}{211}{50.00}{141.00}{212}
\emline{48.00}{141.00}{213}{46.00}{141.00}{214}
\emline{44.00}{141.00}{215}{42.00}{141.00}{216}
\emline{40.00}{141.00}{217}{38.00}{141.00}{218}
\emline{36.00}{141.00}{219}{34.00}{141.00}{220}
\emline{32.00}{141.00}{221}{30.00}{141.00}{222}
\emline{28.00}{141.00}{223}{26.00}{141.00}{224}
\emline{24.00}{141.00}{225}{22.00}{141.00}{226}
\emline{20.00}{141.00}{227}{18.00}{141.00}{228}
\emline{16.00}{141.00}{229}{14.00}{141.00}{230}
\emline{14.00}{139.00}{231}{14.00}{137.00}{232}
\emline{14.00}{135.00}{233}{14.00}{133.00}{234}
\emline{14.00}{131.00}{235}{14.00}{129.00}{236}
\emline{14.00}{127.00}{237}{14.00}{125.00}{238}
\emline{14.00}{123.00}{239}{14.00}{121.00}{240}
\emline{14.00}{119.00}{241}{14.00}{117.00}{242}
\emline{14.00}{115.00}{243}{14.00}{113.00}{244}
\emline{14.00}{111.00}{245}{14.00}{109.00}{246}
\emline{14.00}{107.00}{247}{14.00}{105.00}{248}
\emline{14.00}{103.00}{249}{14.00}{101.00}{250}
\put(35.50,119.00){\makebox(0,0)[cc]{$_\infty$}}
\put(43.50,125.00){\makebox(0,0)[cc]{$_c$}}
\put(47.50,130.00){\makebox(0,0)[cc]{$_4$}}
\put(38.50,130.00){\makebox(0,0)[cc]{$_9$}}
\put(43.50,130.00){\makebox(0,0)[cc]{$_b$}}
\put(43.50,134.00){\makebox(0,0)[cc]{$_3$}}
\put(24.50,115.00){\makebox(0,0)[cc]{$_8$}}
\put(24.50,111.00){\makebox(0,0)[cc]{$_6$}}
\put(20.50,111.00){\makebox(0,0)[cc]{$_5$}}
\put(28.50,111.00){\makebox(0,0)[cc]{$_e$}}
\put(24.50,107.00){\makebox(0,0)[cc]{$_d$}}
\put(39.50,115.00){\makebox(0,0)[cc]{$_1$}}
\put(43.50,110.00){\makebox(0,0)[cc]{$_a$}}
\put(49.50,105.00){\makebox(0,0)[cc]{$_f$}}
\put(18.50,136.00){\makebox(0,0)[cc]{$_2$}}
\put(24.50,130.00){\makebox(0,0)[cc]{$_7$}}
\put(30.50,124.00){\makebox(0,0)[cc]{$_g$}}
\put(55.50,143.00){\makebox(0,0)[cc]{$_0$}}
\emline{137.00}{104.00}{251}{152.00}{104.00}{252}
\emline{152.00}{104.00}{253}{152.00}{119.00}{254}
\emline{152.00}{119.00}{255}{137.00}{119.00}{256}
\emline{137.00}{119.00}{257}{137.00}{104.00}{258}
\emline{156.00}{123.00}{259}{171.00}{123.00}{260}
\emline{171.00}{123.00}{261}{171.00}{138.00}{262}
\emline{171.00}{138.00}{263}{156.00}{138.00}{264}
\emline{156.00}{138.00}{265}{156.00}{123.00}{266}
\emline{137.00}{119.00}{267}{156.00}{138.00}{268}
\emline{152.00}{119.00}{269}{156.00}{123.00}{270}
\emline{152.00}{104.00}{271}{171.00}{123.00}{272}
\emline{171.00}{138.00}{273}{174.00}{141.00}{274}
\emline{137.00}{104.00}{275}{134.00}{101.00}{276}
\emline{136.00}{101.00}{277}{138.00}{101.00}{278}
\emline{140.00}{101.00}{279}{142.00}{101.00}{280}
\emline{144.00}{101.00}{281}{146.00}{101.00}{282}
\emline{148.00}{101.00}{283}{150.00}{101.00}{284}
\emline{152.00}{101.00}{285}{154.00}{101.00}{286}
\emline{156.00}{101.00}{287}{158.00}{101.00}{288}
\emline{160.00}{101.00}{289}{162.00}{101.00}{290}
\emline{164.00}{101.00}{291}{166.00}{101.00}{292}
\emline{168.00}{101.00}{293}{170.00}{101.00}{294}
\emline{172.00}{101.00}{295}{174.00}{101.00}{296}
\emline{174.00}{103.00}{297}{174.00}{105.00}{298}
\emline{174.00}{107.00}{299}{174.00}{109.00}{300}
\emline{174.00}{111.00}{301}{174.00}{113.00}{302}
\emline{174.00}{115.00}{303}{174.00}{117.00}{304}
\emline{174.00}{119.00}{305}{174.00}{121.00}{306}
\emline{174.00}{123.00}{307}{174.00}{125.00}{308}
\emline{174.00}{127.00}{309}{174.00}{129.00}{310}
\emline{174.00}{131.00}{311}{174.00}{133.00}{312}
\emline{174.00}{135.00}{313}{174.00}{137.00}{314}
\emline{174.00}{139.00}{315}{174.00}{141.00}{316}
\emline{172.00}{141.00}{317}{170.00}{141.00}{318}
\emline{168.00}{141.00}{319}{166.00}{141.00}{320}
\emline{164.00}{141.00}{321}{162.00}{141.00}{322}
\emline{160.00}{141.00}{323}{158.00}{141.00}{324}
\emline{156.00}{141.00}{325}{154.00}{141.00}{326}
\emline{152.00}{141.00}{327}{150.00}{141.00}{328}
\emline{148.00}{141.00}{329}{146.00}{141.00}{330}
\emline{144.00}{141.00}{331}{142.00}{141.00}{332}
\emline{140.00}{141.00}{333}{138.00}{141.00}{334}
\emline{136.00}{141.00}{335}{134.00}{141.00}{336}
\emline{134.00}{139.00}{337}{134.00}{137.00}{338}
\emline{134.00}{135.00}{339}{134.00}{133.00}{340}
\emline{134.00}{131.00}{341}{134.00}{129.00}{342}
\emline{134.00}{127.00}{343}{134.00}{125.00}{344}
\emline{134.00}{123.00}{345}{134.00}{121.00}{346}
\emline{134.00}{119.00}{347}{134.00}{117.00}{348}
\emline{134.00}{115.00}{349}{134.00}{113.00}{350}
\emline{134.00}{111.00}{351}{134.00}{109.00}{352}
\emline{134.00}{107.00}{353}{134.00}{105.00}{354}
\emline{134.00}{103.00}{355}{134.00}{101.00}{356}
\put(155.50,119.00){\makebox(0,0)[cc]{$_\infty$}}
\put(163.50,125.00){\makebox(0,0)[cc]{$_3$}}
\put(167.50,130.00){\makebox(0,0)[cc]{$_1$}}
\put(158.50,130.00){\makebox(0,0)[cc]{$_f$}}
\put(163.50,130.00){\makebox(0,0)[cc]{$_7$}}
\put(163.50,134.00){\makebox(0,0)[cc]{$_5$}}
\put(144.50,115.00){\makebox(0,0)[cc]{$_e$}}
\put(144.50,111.00){\makebox(0,0)[cc]{$_a$}}
\put(140.50,111.00){\makebox(0,0)[cc]{$_g$}}
\put(148.50,111.00){\makebox(0,0)[cc]{$_2$}}
\put(144.50,107.00){\makebox(0,0)[cc]{$_c$}}
\put(159.50,115.00){\makebox(0,0)[cc]{$_d$}}
\put(163.50,110.00){\makebox(0,0)[cc]{$_b$}}
\put(169.50,105.00){\makebox(0,0)[cc]{$_8$}}
\put(138.50,136.00){\makebox(0,0)[cc]{$_9$}}
\put(144.50,130.00){\makebox(0,0)[cc]{$_6$}}
\put(150.50,124.00){\makebox(0,0)[cc]{$_4$}}
\put(175.50,143.00){\makebox(0,0)[cc]{$_0$}}
\put(61.50,136.00){\makebox(0,0)[cc]{$\mathbf{_{a^0}}$}}
\put(181.50,136.00){\makebox(0,0)[cc]{$\mathbf{_{c^0}}$}}
\end{picture}
\caption{Label and vertex-tetrahedron representations of
$a^0,b^0,c^0$ in ${\mathcal Q}_3$}
\end{figure}

\begin{figure}[htp]
\hspace*{1.2cm}
\unitlength=0.45mm \special{em:linewidth 0.4pt}
\linethickness{0.4pt}
\begin{picture}(182.00,143.00)
\put(34.00,71.00){\makebox(0,0)[cc]{$_8$}}
\put(34.00,77.00){\makebox(0,0)[cc]{$_b$}}
\put(34.00,65.00){\makebox(0,0)[cc]{$_d$}}
\put(41.00,71.00){\makebox(0,0)[cc]{$_c$}}
\put(28.00,71.00){\makebox(0,0)[cc]{$_\infty$}}
\put(49.00,71.00){\makebox(0,0)[cc]{$_1$}}
\put(50.00,83.00){\makebox(0,0)[cc]{$_f$}}
\put(34.00,56.00){\makebox(0,0)[cc]{$_e$}}
\put(19.00,59.00){\makebox(0,0)[cc]{$_6$}}
\put(19.00,71.00){\makebox(0,0)[cc]{$_9$}}
\put(23.00,86.00){\makebox(0,0)[cc]{$_4$}}
\put(34.00,86.00){\makebox(0,0)[cc]{$_3$}}
\put(57.00,71.00){\makebox(0,0)[cc]{$_0$}}
\put(62.00,71.00){\makebox(0,0)[cc]{$_g$}}
\put(62.00,86.00){\makebox(0,0)[cc]{$\mathbf{_{d^0}}$}}
\put(46.00,56.00){\makebox(0,0)[cc]{$_a$}}
\put(34.00,48.00){\makebox(0,0)[cc]{$_2$}}
\put(34.00,95.00){\makebox(0,0)[cc]{$_7$}}
\put(17.00,77.00){\makebox(0,0)[cc]{$_5$}}
\put(94.00,71.00){\makebox(0,0)[cc]{$_\infty$}}
\put(94.00,77.00){\makebox(0,0)[cc]{$_d$}}
\put(94.00,65.00){\makebox(0,0)[cc]{$_4$}}
\put(101.00,71.00){\makebox(0,0)[cc]{$_1$}}
\put(88.00,71.00){\makebox(0,0)[cc]{$_g$}}
\put(109.00,71.00){\makebox(0,0)[cc]{$_a$}}
\put(111.00,83.00){\makebox(0,0)[cc]{$_5$}}
\put(94.00,56.00){\makebox(0,0)[cc]{$_6$}}
\put(79.00,59.00){\makebox(0,0)[cc]{$_c$}}
\put(79.00,71.00){\makebox(0,0)[cc]{$_7$}}
\put(83.00,86.00){\makebox(0,0)[cc]{$_e$}}
\put(94.00,86.00){\makebox(0,0)[cc]{$_b$}}
\put(117.00,71.00){\makebox(0,0)[cc]{$_f$}}
\put(122.00,71.00){\makebox(0,0)[cc]{$_0$}}
\put(122.00,86.00){\makebox(0,0)[cc]{$\mathbf{_{e^0}}$}}
\put(106.00,56.00){\makebox(0,0)[cc]{$_3$}}
\put(94.00,48.00){\makebox(0,0)[cc]{$_9$}}
\put(94.00,95.00){\makebox(0,0)[cc]{$_8$}}
\put(77.00,77.00){\makebox(0,0)[cc]{$_2$}}
\put(34.00,27.00){\makebox(0,0)[cc]{$_{D_f}$}}
\put(34.00,14.00){\makebox(0,0)[cc]{$_{A_2}$}}
\put(40.00,21.00){\makebox(0,0)[cc]{$_{C_f}$}}
\put(29.00,21.00){\makebox(0,0)[cc]{$_{B_0}$}}
\put(43.00,36.00){\makebox(0,0)[cc]{$_{E_2}$}}
\put(25.00,5.00){\makebox(0,0)[cc]{$_{D_2}$}}
\put(25.00,36.00){\makebox(0,0)[cc]{$_{A_f}$}}
\put(50.00,21.00){\makebox(0,0)[cc]{$_{E_0}$}}
\put(44.00,5.00){\makebox(0,0)[cc]{$_{A_g}$}}
\put(35.00,4.00){\makebox(0,0)[cc]{$_{E_f}$}}
\put(35.00,37.00){\makebox(0,0)[cc]{$_{A_1}$}}
\put(18.00,21.00){\makebox(0,0)[cc]{$_{C_2}$}}
\put(93.00,27.00){\makebox(0,0)[cc]{$_{B_b}$}}
\put(93.00,14.00){\makebox(0,0)[cc]{$_{B_6}$}}
\put(100.00,21.00){\makebox(0,0)[cc]{$_{B_a}$}}
\put(89.00,21.00){\makebox(0,0)[cc]{$_{B_7}$}}
\put(103.00,36.00){\makebox(0,0)[cc]{$_{C_2}$}}
\put(85.00,5.00){\makebox(0,0)[cc]{$_{C_f}$}}
\put(85.00,36.00){\makebox(0,0)[cc]{$_{A_9}$}}
\put(110.00,21.00){\makebox(0,0)[cc]{$_{C_7}$}}
\put(103.00,5.00){\makebox(0,0)[cc]{$_{A_8}$}}
\put(95.00,4.00){\makebox(0,0)[cc]{$_{A_b}$}}
\put(95.00,37.00){\makebox(0,0)[cc]{$_{A_6}$}}
\put(78.00,21.00){\makebox(0,0)[cc]{$_{C_a}$}}
\put(60.00,41.00){\makebox(0,0)[cc]{$_{C^1}$}}
\put(60.00,1.00){\makebox(0,0)[cc]{$_{C^g}$}}
\put(9.00,41.00){\makebox(0,0)[cc]{$_{E^8}$}}
\put(9.00,1.00){\makebox(0,0)[cc]{$_{A^2}$}}
\put(19.00,29.00){\makebox(0,0)[cc]{$_{D^f}$}}
\put(19.00,11.00){\makebox(0,0)[cc]{$_{D^2}$}}
\put(50.00,31.00){\makebox(0,0)[cc]{$_{A^f}$}}
\put(50.00,11.00){\makebox(0,0)[cc]{$_{E^9}$}}
\put(120.00,41.00){\makebox(0,0)[cc]{$_{E^d}$}}
\put(120.00,1.00){\makebox(0,0)[cc]{$_{A^b}$}}
\put(69.00,41.00){\makebox(0,0)[cc]{$_{E^g}$}}
\put(69.00,1.00){\makebox(0,0)[cc]{$_{E^4}$}}
\put(79.00,29.00){\makebox(0,0)[cc]{$_{D^9}$}}
\put(79.00,11.00){\makebox(0,0)[cc]{$_{F^f}$}}
\put(110.00,31.00){\makebox(0,0)[cc]{$_{F^2}$}}
\put(110.00,11.00){\makebox(0,0)[cc]{$_{D^a}$}}
\put(129.00,41.00){\makebox(0,0)[cc]{$_{E^f}$}}
\put(129.00,1.00){\makebox(0,0)[cc]{$_{C^8}$}}
\emline{14.00}{51.00}{1}{14.00}{91.00}{2}
\emline{14.00}{91.00}{3}{54.00}{91.00}{4}
\emline{54.00}{91.00}{5}{44.00}{81.00}{6}
\emline{44.00}{81.00}{7}{24.00}{81.00}{8}
\emline{44.00}{81.00}{9}{44.00}{61.00}{10}
\emline{14.00}{91.00}{11}{24.00}{81.00}{12}
\emline{14.00}{51.00}{13}{54.00}{51.00}{14}
\emline{54.00}{51.00}{15}{54.00}{91.00}{16}
\emline{24.00}{81.00}{17}{24.00}{61.00}{18}
\emline{44.00}{61.00}{19}{24.00}{61.00}{20}
\emline{44.00}{61.00}{21}{54.00}{51.00}{22}
\emline{24.00}{61.00}{23}{14.00}{51.00}{24}
\emline{14.00}{1.00}{25}{14.00}{41.00}{26}
\emline{14.00}{41.00}{27}{54.00}{41.00}{28}
\emline{54.00}{41.00}{29}{44.00}{31.00}{30}
\emline{44.00}{31.00}{31}{24.00}{31.00}{32}
\emline{44.00}{31.00}{33}{44.00}{11.00}{34}
\emline{14.00}{41.00}{35}{24.00}{31.00}{36}
\emline{14.00}{1.00}{37}{54.00}{1.00}{38}
\emline{54.00}{1.00}{39}{54.00}{41.00}{40}
\emline{24.00}{31.00}{41}{24.00}{11.00}{42}
\emline{44.00}{11.00}{43}{24.00}{11.00}{44}
\emline{44.00}{11.00}{45}{54.00}{1.00}{46}
\emline{24.00}{11.00}{47}{14.00}{1.00}{48}
\emline{74.00}{51.00}{49}{74.00}{91.00}{50}
\emline{74.00}{91.00}{51}{114.00}{91.00}{52}
\emline{114.00}{91.00}{53}{104.00}{81.00}{54}
\emline{104.00}{81.00}{55}{84.00}{81.00}{56}
\emline{104.00}{81.00}{57}{104.00}{61.00}{58}
\emline{74.00}{91.00}{59}{84.00}{81.00}{60}
\emline{74.00}{51.00}{61}{114.00}{51.00}{62}
\emline{114.00}{51.00}{63}{114.00}{91.00}{64}
\emline{84.00}{81.00}{65}{84.00}{61.00}{66}
\emline{104.00}{61.00}{67}{84.00}{61.00}{68}
\emline{104.00}{61.00}{69}{114.00}{51.00}{70}
\emline{84.00}{61.00}{71}{74.00}{51.00}{72}
\emline{74.00}{1.00}{73}{74.00}{41.00}{74}
\emline{74.00}{41.00}{75}{114.00}{41.00}{76}
\emline{114.00}{41.00}{77}{104.00}{31.00}{78}
\emline{104.00}{31.00}{79}{84.00}{31.00}{80}
\emline{104.00}{31.00}{81}{104.00}{11.00}{82}
\emline{74.00}{41.00}{83}{84.00}{31.00}{84}
\emline{74.00}{1.00}{85}{114.00}{1.00}{86}
\emline{114.00}{1.00}{87}{114.00}{41.00}{88}
\emline{84.00}{31.00}{89}{84.00}{11.00}{90}
\emline{104.00}{11.00}{91}{84.00}{11.00}{92}
\emline{104.00}{11.00}{93}{114.00}{1.00}{94}
\emline{84.00}{11.00}{95}{74.00}{1.00}{96}
\emline{17.50}{104.00}{97}{32.50}{104.00}{98}
\emline{32.50}{104.00}{99}{32.50}{119.00}{100}
\emline{32.50}{119.00}{101}{17.50}{119.00}{102}
\emline{17.50}{119.00}{103}{17.50}{104.00}{104}
\emline{36.50}{123.00}{105}{51.50}{123.00}{106}
\emline{51.50}{123.00}{107}{51.50}{138.00}{108}
\emline{51.50}{138.00}{109}{36.50}{138.00}{110}
\emline{36.50}{138.00}{111}{36.50}{123.00}{112}
\emline{17.50}{119.00}{113}{36.50}{138.00}{114}
\emline{32.50}{119.00}{115}{36.50}{123.00}{116}
\emline{32.50}{104.00}{117}{51.50}{123.00}{118}
\emline{51.50}{138.00}{119}{54.50}{141.00}{120}
\emline{17.50}{104.00}{121}{14.50}{101.00}{122}
\emline{16.50}{101.00}{123}{18.50}{101.00}{124}
\emline{20.50}{101.00}{125}{22.50}{101.00}{126}
\emline{24.50}{101.00}{127}{26.50}{101.00}{128}
\emline{28.50}{101.00}{129}{30.50}{101.00}{130}
\emline{32.50}{101.00}{131}{34.50}{101.00}{132}
\emline{36.50}{101.00}{133}{38.50}{101.00}{134}
\emline{40.50}{101.00}{135}{42.50}{101.00}{136}
\emline{44.50}{101.00}{137}{46.50}{101.00}{138}
\emline{48.50}{101.00}{139}{50.50}{101.00}{140}
\emline{52.50}{101.00}{141}{54.50}{101.00}{142}
\emline{54.50}{103.00}{143}{54.50}{105.00}{144}
\emline{54.50}{107.00}{145}{54.50}{109.00}{146}
\emline{54.50}{111.00}{147}{54.50}{113.00}{148}
\emline{54.50}{115.00}{149}{54.50}{117.00}{150}
\emline{54.50}{119.00}{151}{54.50}{121.00}{152}
\emline{54.50}{123.00}{153}{54.50}{125.00}{154}
\emline{54.50}{127.00}{155}{54.50}{129.00}{156}
\emline{54.50}{131.00}{157}{54.50}{133.00}{158}
\emline{54.50}{135.00}{159}{54.50}{137.00}{160}
\emline{54.50}{139.00}{161}{54.50}{141.00}{162}
\emline{52.50}{141.00}{163}{50.50}{141.00}{164}
\emline{48.50}{141.00}{165}{46.50}{141.00}{166}
\emline{44.50}{141.00}{167}{42.50}{141.00}{168}
\emline{40.50}{141.00}{169}{38.50}{141.00}{170}
\emline{36.50}{141.00}{171}{34.50}{141.00}{172}
\emline{32.50}{141.00}{173}{30.50}{141.00}{174}
\emline{28.50}{141.00}{175}{26.50}{141.00}{176}
\emline{24.50}{141.00}{177}{22.50}{141.00}{178}
\emline{20.50}{141.00}{179}{18.50}{141.00}{180}
\emline{16.50}{141.00}{181}{14.50}{141.00}{182}
\emline{14.50}{139.00}{183}{14.50}{137.00}{184}
\emline{14.50}{135.00}{185}{14.50}{133.00}{186}
\emline{14.50}{131.00}{187}{14.50}{129.00}{188}
\emline{14.50}{127.00}{189}{14.50}{125.00}{190}
\emline{14.50}{123.00}{191}{14.50}{121.00}{192}
\emline{14.50}{119.00}{193}{14.50}{117.00}{194}
\emline{14.50}{115.00}{195}{14.50}{113.00}{196}
\emline{14.50}{111.00}{197}{14.50}{109.00}{198}
\emline{14.50}{107.00}{199}{14.50}{105.00}{200}
\emline{14.50}{103.00}{201}{14.50}{101.00}{202}
\put(36.00,119.00){\makebox(0,0)[cc]{$_\infty$}}
\put(44.00,125.00){\makebox(0,0)[cc]{$_b$}}
\put(48.00,130.00){\makebox(0,0)[cc]{$_f$}}
\put(39.00,130.00){\makebox(0,0)[cc]{$_4$}}
\put(44.00,130.00){\makebox(0,0)[cc]{$_3$}}
\put(44.00,134.00){\makebox(0,0)[cc]{$_7$}}
\put(25.00,115.00){\makebox(0,0)[cc]{$_6$}}
\put(25.00,111.00){\makebox(0,0)[cc]{$_e$}}
\put(21.00,111.00){\makebox(0,0)[cc]{$_2$}}
\put(29.00,111.00){\makebox(0,0)[cc]{$_d$}}
\put(25.00,107.00){\makebox(0,0)[cc]{$_a$}}
\put(40.00,115.00){\makebox(0,0)[cc]{$_8$}}
\put(44.00,110.00){\makebox(0,0)[cc]{$_c$}}
\put(50.00,105.00){\makebox(0,0)[cc]{$_1$}}
\put(19.00,136.00){\makebox(0,0)[cc]{$_g$}}
\put(25.00,130.00){\makebox(0,0)[cc]{$_5$}}
\put(31.00,124.00){\makebox(0,0)[cc]{$_9$}}
\put(56.00,143.00){\makebox(0,0)[cc]{$_0$}}
\put(62.00,136.00){\makebox(0,0)[cc]{$\mathbf{_{d^0}}$}}
\put(154.00,71.00){\makebox(0,0)[cc]{$_f$}}
\put(154.00,77.00){\makebox(0,0)[cc]{$_a$}}
\put(154.00,65.00){\makebox(0,0)[cc]{$_1$}}
\put(161.00,71.00){\makebox(0,0)[cc]{$_e$}}
\put(148.00,71.00){\makebox(0,0)[cc]{$_\infty$}}
\put(169.00,71.00){\makebox(0,0)[cc]{$_4$}}
\put(171.00,83.00){\makebox(0,0)[cc]{$_9$}}
\put(154.00,56.00){\makebox(0,0)[cc]{$_5$}}
\put(139.00,59.00){\makebox(0,0)[cc]{$_7$}}
\put(139.00,71.00){\makebox(0,0)[cc]{$_2$}}
\put(143.00,86.00){\makebox(0,0)[cc]{$_g$}}
\put(154.00,86.00){\makebox(0,0)[cc]{$_c$}}
\put(177.00,71.00){\makebox(0,0)[cc]{$_0$}}
\put(182.00,71.00){\makebox(0,0)[cc]{$_d$}}
\put(182.00,86.00){\makebox(0,0)[cc]{$\mathbf{_{f^0}}$}}
\put(166.00,56.00){\makebox(0,0)[cc]{$_6$}}
\put(154.00,48.00){\makebox(0,0)[cc]{$_8$}}
\put(154.00,95.00){\makebox(0,0)[cc]{$_b$}}
\put(137.00,77.00){\makebox(0,0)[cc]{$_3$}}
\put(153.00,27.00){\makebox(0,0)[cc]{$_{F_9}$}}
\put(153.00,14.00){\makebox(0,0)[cc]{$_{C_8}$}}
\put(160.00,21.00){\makebox(0,0)[cc]{$_{A_9}$}}
\put(149.00,21.00){\makebox(0,0)[cc]{$_{B_0}$}}
\put(163.00,36.00){\makebox(0,0)[cc]{$_{E_8}$}}
\put(145.00,5.00){\makebox(0,0)[cc]{$_{F_8}$}}
\put(145.00,36.00){\makebox(0,0)[cc]{$_{C_9}$}}
\put(170.00,21.00){\makebox(0,0)[cc]{$_{E_0}$}}
\put(163.00,5.00){\makebox(0,0)[cc]{$_{C_d}$}}
\put(155.00,4.00){\makebox(0,0)[cc]{$_{E_9}$}}
\put(155.00,37.00){\makebox(0,0)[cc]{$_{C_4}$}}
\put(138.00,21.00){\makebox(0,0)[cc]{$_{A_8}$}}
\put(180.00,41.00){\makebox(0,0)[cc]{$_{A^4}$}}
\put(180.00,1.00){\makebox(0,0)[cc]{$_{A^d}$}}
\put(139.00,30.00){\makebox(0,0)[cc]{$_{F^9}$}}
\put(139.00,11.00){\makebox(0,0)[cc]{$_{F^8}$}}
\put(170.00,31.00){\makebox(0,0)[cc]{$_{C^9}$}}
\put(170.00,11.00){\makebox(0,0)[cc]{$_{E^2}$}}
\emline{134.00}{51.00}{203}{134.00}{91.00}{204}
\emline{134.00}{91.00}{205}{174.00}{91.00}{206}
\emline{174.00}{91.00}{207}{164.00}{81.00}{208}
\emline{164.00}{81.00}{209}{144.00}{81.00}{210}
\emline{164.00}{81.00}{211}{164.00}{61.00}{212}
\emline{134.00}{91.00}{213}{144.00}{81.00}{214}
\emline{134.00}{51.00}{215}{174.00}{51.00}{216}
\emline{174.00}{51.00}{217}{174.00}{91.00}{218}
\emline{144.00}{81.00}{219}{144.00}{61.00}{220}
\emline{164.00}{61.00}{221}{144.00}{61.00}{222}
\emline{164.00}{61.00}{223}{174.00}{51.00}{224}
\emline{144.00}{61.00}{225}{134.00}{51.00}{226}
\emline{134.00}{1.00}{227}{134.00}{41.00}{228}
\emline{134.00}{41.00}{229}{174.00}{41.00}{230}
\emline{174.00}{41.00}{231}{164.00}{31.00}{232}
\emline{164.00}{31.00}{233}{144.00}{31.00}{234}
\emline{164.00}{31.00}{235}{164.00}{11.00}{236}
\emline{134.00}{41.00}{237}{144.00}{31.00}{238}
\emline{134.00}{1.00}{239}{174.00}{1.00}{240}
\emline{174.00}{1.00}{241}{174.00}{41.00}{242}
\emline{144.00}{31.00}{243}{144.00}{11.00}{244}
\emline{164.00}{11.00}{245}{144.00}{11.00}{246}
\emline{164.00}{11.00}{247}{174.00}{1.00}{248}
\emline{144.00}{11.00}{249}{134.00}{1.00}{250}
\emline{137.50}{104.00}{251}{152.50}{104.00}{252}
\emline{152.50}{104.00}{253}{152.50}{119.00}{254}
\emline{152.50}{119.00}{255}{137.50}{119.00}{256}
\emline{137.50}{119.00}{257}{137.50}{104.00}{258}
\emline{156.50}{123.00}{259}{171.50}{123.00}{260}
\emline{171.50}{123.00}{261}{171.50}{138.00}{262}
\emline{171.50}{138.00}{263}{156.50}{138.00}{264}
\emline{156.50}{138.00}{265}{156.50}{123.00}{266}
\emline{137.50}{119.00}{267}{156.50}{138.00}{268}
\emline{152.50}{119.00}{269}{156.50}{123.00}{270}
\emline{152.50}{104.00}{271}{171.50}{123.00}{272}
\emline{171.50}{138.00}{273}{174.50}{141.00}{274}
\emline{137.50}{104.00}{275}{134.50}{101.00}{276}
\emline{136.50}{101.00}{277}{138.50}{101.00}{278}
\emline{140.50}{101.00}{279}{142.50}{101.00}{280}
\emline{144.50}{101.00}{281}{146.50}{101.00}{282}
\emline{148.50}{101.00}{283}{150.50}{101.00}{284}
\emline{152.50}{101.00}{285}{154.50}{101.00}{286}
\emline{156.50}{101.00}{287}{158.50}{101.00}{288}
\emline{160.50}{101.00}{289}{162.50}{101.00}{290}
\emline{164.50}{101.00}{291}{166.50}{101.00}{292}
\emline{168.50}{101.00}{293}{170.50}{101.00}{294}
\emline{172.50}{101.00}{295}{174.50}{101.00}{296}
\emline{174.50}{103.00}{297}{174.50}{105.00}{298}
\emline{174.50}{107.00}{299}{174.50}{109.00}{300}
\emline{174.50}{111.00}{301}{174.50}{113.00}{302}
\emline{174.50}{115.00}{303}{174.50}{117.00}{304}
\emline{174.50}{119.00}{305}{174.50}{121.00}{306}
\emline{174.50}{123.00}{307}{174.50}{125.00}{308}
\emline{174.50}{127.00}{309}{174.50}{129.00}{310}
\emline{174.50}{131.00}{311}{174.50}{133.00}{312}
\emline{174.50}{135.00}{313}{174.50}{137.00}{314}
\emline{174.50}{139.00}{315}{174.50}{141.00}{316}
\emline{172.50}{141.00}{317}{170.50}{141.00}{318}
\emline{168.50}{141.00}{319}{166.50}{141.00}{320}
\emline{164.50}{141.00}{321}{162.50}{141.00}{322}
\emline{160.50}{141.00}{323}{158.50}{141.00}{324}
\emline{156.50}{141.00}{325}{154.50}{141.00}{326}
\emline{152.50}{141.00}{327}{150.50}{141.00}{328}
\emline{148.50}{141.00}{329}{146.50}{141.00}{330}
\emline{144.50}{141.00}{331}{142.50}{141.00}{332}
\emline{140.50}{141.00}{333}{138.50}{141.00}{334}
\emline{136.50}{141.00}{335}{134.50}{141.00}{336}
\emline{134.50}{139.00}{337}{134.50}{137.00}{338}
\emline{134.50}{135.00}{339}{134.50}{133.00}{340}
\emline{134.50}{131.00}{341}{134.50}{129.00}{342}
\emline{134.50}{127.00}{343}{134.50}{125.00}{344}
\emline{134.50}{123.00}{345}{134.50}{121.00}{346}
\emline{134.50}{119.00}{347}{134.50}{117.00}{348}
\emline{134.50}{115.00}{349}{134.50}{113.00}{350}
\emline{134.50}{111.00}{351}{134.50}{109.00}{352}
\emline{134.50}{107.00}{353}{134.50}{105.00}{354}
\put(156.00,119.00){\makebox(0,0)[cc]{$_\infty$}}
\put(164.00,125.00){\makebox(0,0)[cc]{$_a$}}
\put(168.00,130.00){\makebox(0,0)[cc]{$_9$}}
\put(159.00,130.00){\makebox(0,0)[cc]{$_g$}}
\put(164.00,130.00){\makebox(0,0)[cc]{$_c$}}
\put(164.00,134.00){\makebox(0,0)[cc]{$_b$}}
\put(145.00,115.00){\makebox(0,0)[cc]{$_7$}}
\put(145.00,111.00){\makebox(0,0)[cc]{$_5$}}
\put(141.00,111.00){\makebox(0,0)[cc]{$_8$}}
\put(149.00,111.00){\makebox(0,0)[cc]{$_1$}}
\put(145.00,107.00){\makebox(0,0)[cc]{$_6$}}
\put(160.00,115.00){\makebox(0,0)[cc]{$_f$}}
\put(164.00,110.00){\makebox(0,0)[cc]{$_e$}}
\put(170.00,105.00){\makebox(0,0)[cc]{$_4$}}
\put(139.00,136.00){\makebox(0,0)[cc]{$_d$}}
\put(145.00,130.00){\makebox(0,0)[cc]{$_3$}}
\put(151.00,124.00){\makebox(0,0)[cc]{$_2$}}
\put(176.00,143.00){\makebox(0,0)[cc]{$_0$}}
\put(182.00,136.00){\makebox(0,0)[cc]{$\mathbf{_{f^0}}$}}
\emline{134.50}{103.00}{355}{134.50}{101.00}{356}
\end{picture}
\caption{Label and vertex-tetrahedron representations of
$d^0,e^0,f^0$ in ${\mathcal Q}_3$}
\end{figure}
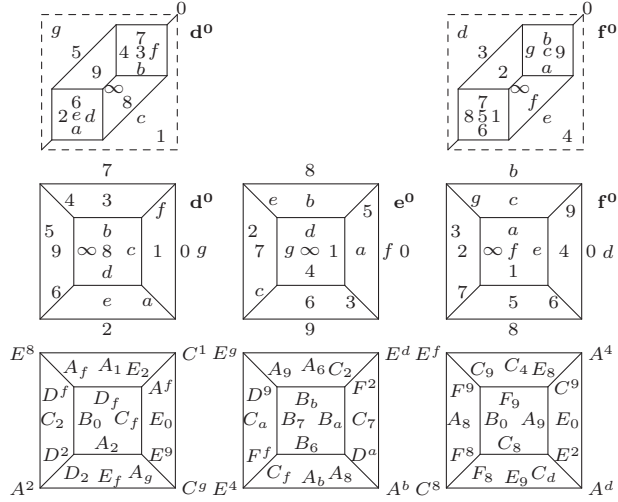

\noindent The labels of the 12 vertices and 6 4-holes of each of
$\sigma^0=a^0,\ldots,f^0$ are depicted again on the middle thirds of
Figures 5 and 6, this time on a copy ${\mathcal Q}_3$ of the 3-cube
$Q_3$ from which a corresponding copy of $L(Q_3)$ in $\mathcal Y$ is
obtained with its vertices taken as the middle points of the edges
of ${\mathcal Q}_3$, tracing an edge between two such vertices
whenever the edges they represent have a vertex in common in
${\mathcal Q}_3$, with the convention that labels of vertices and
4-holes of $\sigma^0$ label now respectively the corresponding edges
and faces of ${\mathcal Q}_3$. (On the bottom thirds those edges are
labeled by the corresponding vertices of $\mathcal S$ and their
vertices by the corresponding containing copies of $K_4$; on the
upper thirds, 4 different cutouts of ${\mathcal Q}_3$ are depicted
to show involution symmetry around edges labeled $\infty$, where
${\mathcal Q}_3$ is regained by identifying the upper and left sides
and the lower and right sides via 90$^o$ rotations at the upper-left
and lower-right corners). Opposite faces in such $\sigma^j$
determine pairs of points of $PG(1,17)$, a total of 3 such pairs
leading to a unique sextet which is not a vertex of $\mathcal S$ but
uniformly 3 times a vertex of $\mathcal S$. For example, these 3
pairs for $\sigma^0=a^0$ form the sextet
$\{12,6b,fg\}=3\times\{6c,2f,5b\}=A_0$, mod 17. By denoting
$a^0=\{12,6b,fg\}$ and so on for the 101 remaining copies of
$L(Q_3)$ in $PG(1,17)$, we obtain a self-dual configuration that
uses again the duality $\phi$ of Section 4, this time with points
and lines taken as the vertices and copies of $L(Q_3)$ in $\mathcal
S$, as claimed in Theorem 6.1{\bf(8)} below.

\begin{figure}[htp]
\hspace*{.1mm}
\unitlength=0.50mm \special{em:linewidth 0.4pt}
\linethickness{0.4pt}
\begin{picture}(216.00,49.00)
\put(11.00,15.00){\circle{2.00}} \put(21.00,25.00){\circle{2.00}}
\put(31.00,15.00){\circle{2.00}}
\emline{2.00}{5.00}{1}{20.00}{5.00}{2}
\emline{22.00}{5.00}{3}{40.00}{5.00}{4}
\put(11.00,35.00){\circle{2.00}} \put(31.00,35.00){\circle{2.00}}
\put(31.00,35.00){\circle{2.00}}
\emline{2.00}{45.00}{5}{20.00}{45.00}{6}
\emline{22.00}{45.00}{7}{40.00}{45.00}{8}
\put(51.00,15.00){\circle{2.00}} \put(61.00,25.00){\circle{2.00}}
\put(71.00,15.00){\circle{2.00}}
\emline{42.00}{5.00}{9}{60.00}{5.00}{10}
\emline{62.00}{5.00}{11}{80.00}{5.00}{12}
\put(51.00,35.00){\circle{2.00}} \put(71.00,35.00){\circle{2.00}}
\put(71.00,35.00){\circle{2.00}}
\emline{42.00}{45.00}{13}{60.00}{45.00}{14}
\emline{62.00}{45.00}{15}{80.00}{45.00}{16}
\put(91.00,15.00){\circle{2.00}} \put(101.00,25.00){\circle{2.00}}
\put(111.00,15.00){\circle{2.00}}
\emline{82.00}{5.00}{17}{100.00}{5.00}{18}
\emline{102.00}{5.00}{19}{120.00}{5.00}{20}
\put(91.00,35.00){\circle{2.00}} \put(111.00,35.00){\circle{2.00}}
\put(111.00,35.00){\circle{2.00}}
\emline{82.00}{45.00}{21}{100.00}{45.00}{22}
\emline{102.00}{45.00}{23}{120.00}{45.00}{24}
\put(131.00,15.00){\circle{2.00}} \put(141.00,25.00){\circle{2.00}}
\put(151.00,15.00){\circle{2.00}}
\emline{122.00}{5.00}{25}{140.00}{5.00}{26}
\emline{142.00}{5.00}{27}{160.00}{5.00}{28}
\put(131.00,35.00){\circle{2.00}} \put(151.00,35.00){\circle{2.00}}
\emline{122.00}{45.00}{29}{140.00}{45.00}{30}
\emline{142.00}{45.00}{31}{160.00}{45.00}{32}
\put(3.00,49.00){\makebox(0,0)[cc]{$_{C_d^{3g}}$}}
\put(22.00,49.00){\makebox(0,0)[cc]{$_{B_5^7}$}}
\put(42.00,49.00){\makebox(0,0)[cc]{$_{A_3^{89}}$}}
\put(62.00,49.00){\makebox(0,0)[cc]{$_{E_4^0}$}}
\put(82.00,49.00){\makebox(0,0)[cc]{$_{D_0^{d7}}$}}
\put(102.00,49.00){\makebox(0,0)[cc]{$_{F_2^g}$}}
\put(122.00,49.00){\makebox(0,0)[cc]{$_{E_1^{\infty 0}}$}}
\put(142.00,49.00){\makebox(0,0)[cc]{$_{E_9^9}$}}
\put(3.00,1.00){\makebox(0,0)[cc]{$_{C_1^{89}}$}}
\put(21.00,1.00){\makebox(0,0)[cc]{$_{A_5^3}$}}
\put(41.00,1.00){\makebox(0,0)[cc]{$_{C_4^{d7}}$}}
\put(61.00,1.00){\makebox(0,0)[cc]{$_{D_4^8}$}}
\put(81.00,1.00){\makebox(0,0)[cc]{$_{B_2^{\infty 0}}$}}
\put(101.00,1.00){\makebox(0,0)[cc]{$_{F_a^d}$}}
\put(121.00,1.00){\makebox(0,0)[cc]{$_{F_0^{3g}}$}}
\put(141.00,1.00){\makebox(0,0)[cc]{$_{B_9^\infty}$}}
\put(6.00,35.00){\makebox(0,0)[cc]{$_{B_4^c}$}}
\put(26.00,35.00){\makebox(0,0)[cc]{$_{B_1^a}$}}
\put(46.00,35.00){\makebox(0,0)[cc]{$_{E_2^2}$}}
\put(66.00,35.00){\makebox(0,0)[cc]{$_{C_0^e}$}}
\put(86.00,35.00){\makebox(0,0)[cc]{$_{F_9^6}$}}
\put(106.00,35.00){\makebox(0,0)[cc]{$_{A_2^4}$}}
\put(126.00,35.00){\makebox(0,0)[cc]{$_{C_5^f}$}}
\put(146.00,35.00){\makebox(0,0)[cc]{$_{E_0^1}$}}
\put(6.00,15.00){\makebox(0,0)[cc]{$_{A_2^e}$}}
\put(26.00,15.00){\makebox(0,0)[cc]{$_{C_9^6}$}}
\put(46.00,15.00){\makebox(0,0)[cc]{$_{E_0^4}$}}
\put(66.00,15.00){\makebox(0,0)[cc]{$_{A_4^f}$}}
\put(86.00,15.00){\makebox(0,0)[cc]{$_{B_1^1}$}}
\put(106.00,15.00){\makebox(0,0)[cc]{$_{D_2^c}$}}
\put(126.00,15.00){\makebox(0,0)[cc]{$_{C_0^a}$}}
\put(146.00,15.00){\makebox(0,0)[cc]{$_{F_1^2}$}}
\put(21.00,20.00){\makebox(0,0)[cc]{$_{B_0^f}$}}
\put(61.00,20.00){\makebox(0,0)[cc]{$_{A_1^c}$}}
\put(163.00,49.00){\makebox(0,0)[cc]{$_{C_d^{3g}}$}}
\put(161.00,1.00){\makebox(0,0)[cc]{$_{C_1^{89}}$}}
\put(4.00,29.00){\makebox(0,0)[cc]{$_{E^7}$}}
\put(21.00,31.00){\makebox(0,0)[cc]{$_{e^b}$}}
\put(44.00,29.00){\makebox(0,0)[cc]{$_{E^a}$}}
\put(61.00,31.00){\makebox(0,0)[cc]{$_{d^2}$}}
\put(84.00,29.00){\makebox(0,0)[cc]{$_{D^0}$}}
\put(101.00,31.00){\makebox(0,0)[cc]{$_{c^1}$}}
\put(124.00,29.00){\makebox(0,0)[cc]{$_{C^0}$}}
\put(141.00,31.00){\makebox(0,0)[cc]{$_{f^{\,9}}$}}
\put(158.74,30.00){\makebox(0,0)[cc]{$_{E^7}$}}
\put(1.00,45.00){\circle*{2.00}} \put(21.00,45.00){\circle*{2.00}}
\put(41.00,45.00){\circle*{2.00}} \put(61.00,45.00){\circle*{2.00}}
\put(1.00,5.00){\circle*{2.00}} \put(21.00,5.00){\circle*{2.00}}
\put(41.00,5.00){\circle*{2.00}} \put(61.00,5.00){\circle*{2.00}}
\put(81.00,45.00){\circle*{2.00}} \put(101.00,45.00){\circle*{2.00}}
\put(121.00,45.00){\circle*{2.00}}
\put(141.00,45.00){\circle*{2.00}} \put(81.00,5.00){\circle*{2.00}}
\put(101.00,5.00){\circle*{2.00}} \put(121.00,5.00){\circle*{2.00}}
\put(141.00,5.00){\circle*{2.00}} \put(161.00,45.00){\circle*{2.00}}
\put(161.00,5.00){\circle*{2.00}}
\emline{2.00}{44.00}{33}{10.00}{36.00}{34}
\emline{12.00}{36.00}{35}{20.00}{44.00}{36}
\emline{12.00}{34.00}{37}{20.00}{26.00}{38}
\emline{20.00}{24.00}{39}{12.00}{16.00}{40}
\emline{10.00}{14.00}{41}{2.00}{6.00}{42}
\emline{12.00}{14.00}{43}{20.00}{6.00}{44}
\emline{11.00}{34.00}{45}{11.00}{16.00}{46}
\emline{22.00}{44.00}{47}{30.00}{36.00}{48}
\emline{30.00}{34.00}{49}{22.00}{26.00}{50}
\emline{22.00}{24.00}{51}{30.00}{16.00}{52}
\emline{30.00}{14.00}{53}{22.00}{6.00}{54}
\emline{1.00}{44.00}{55}{1.00}{6.00}{56}
\emline{31.00}{34.00}{57}{31.00}{16.00}{58}
\emline{32.00}{36.00}{59}{40.00}{44.00}{60}
\emline{32.00}{14.00}{61}{40.00}{6.00}{62}
\emline{42.00}{44.00}{63}{50.00}{36.00}{64}
\emline{52.00}{36.00}{65}{60.00}{44.00}{66}
\emline{60.00}{24.00}{67}{52.00}{16.00}{68}
\emline{50.00}{14.00}{69}{42.00}{6.00}{70}
\emline{52.00}{14.00}{71}{60.00}{6.00}{72}
\emline{51.00}{34.00}{73}{51.00}{16.00}{74}
\emline{62.00}{44.00}{75}{70.00}{36.00}{76}
\emline{70.00}{34.00}{77}{62.00}{26.00}{78}
\emline{62.00}{24.00}{79}{70.00}{16.00}{80}
\emline{70.00}{14.00}{81}{62.00}{6.00}{82}
\emline{41.00}{44.00}{83}{41.00}{6.00}{84}
\emline{71.00}{34.00}{85}{71.00}{16.00}{86}
\emline{72.00}{36.00}{87}{80.00}{44.00}{88}
\emline{72.00}{14.00}{89}{80.00}{6.00}{90}
\emline{52.00}{34.00}{91}{60.00}{26.00}{92}
\emline{82.00}{44.00}{93}{90.00}{36.00}{94}
\emline{92.00}{36.00}{95}{100.00}{44.00}{96}
\emline{92.00}{34.00}{97}{100.00}{26.00}{98}
\emline{100.00}{24.00}{99}{92.00}{16.00}{100}
\emline{90.00}{14.00}{101}{82.00}{6.00}{102}
\emline{92.00}{14.00}{103}{100.00}{6.00}{104}
\emline{91.00}{34.00}{105}{91.00}{16.00}{106}
\emline{102.00}{44.00}{107}{110.00}{36.00}{108}
\emline{110.00}{34.00}{109}{102.00}{26.00}{110}
\emline{102.00}{24.00}{111}{110.00}{16.00}{112}
\emline{110.00}{14.00}{113}{102.00}{6.00}{114}
\emline{81.00}{44.00}{115}{81.00}{6.00}{116}
\emline{111.00}{34.00}{117}{111.00}{16.00}{118}
\emline{112.00}{36.00}{119}{120.00}{44.00}{120}
\emline{112.00}{14.00}{121}{120.00}{6.00}{122}
\emline{122.00}{44.00}{123}{130.00}{36.00}{124}
\emline{132.00}{36.00}{125}{140.00}{44.00}{126}
\emline{140.00}{24.00}{127}{132.00}{16.00}{128}
\emline{130.00}{14.00}{129}{122.00}{6.00}{130}
\emline{132.00}{14.00}{131}{140.00}{6.00}{132}
\emline{131.00}{34.00}{133}{131.00}{16.00}{134}
\emline{142.00}{44.00}{135}{150.00}{36.00}{136}
\emline{150.00}{34.00}{137}{142.00}{26.00}{138}
\emline{142.00}{24.00}{139}{150.00}{16.00}{140}
\emline{150.00}{14.00}{141}{142.00}{6.00}{142}
\emline{121.00}{44.00}{143}{121.00}{6.00}{144}
\emline{151.00}{34.00}{145}{151.00}{16.00}{146}
\emline{152.00}{36.00}{147}{160.00}{44.00}{148}
\emline{152.00}{14.00}{149}{160.00}{6.00}{150}
\emline{132.00}{34.00}{151}{140.00}{26.00}{152}
\emline{161.00}{44.00}{153}{161.00}{6.00}{154}
\put(19.00,38.00){\makebox(0,0)[cc]{$_\infty$}}
\put(5.00,23.00){\makebox(0,0)[cc]{$_4$}}
\put(19.00,13.00){\makebox(0,0)[cc]{$_0$}}
\put(36.00,23.00){\makebox(0,0)[cc]{$_1$}}
\put(31.00,49.00){\makebox(0,0)[cc]{$_5$}}
\put(31.00,1.00){\makebox(0,0)[cc]{$_b$}}
\put(59.00,38.00){\makebox(0,0)[cc]{$_3$}}
\put(59.00,13.00){\makebox(0,0)[cc]{$_g$}}
\put(76.00,23.00){\makebox(0,0)[cc]{$_a$}}
\put(71.00,49.00){\makebox(0,0)[cc]{$_5$}}
\put(71.00,1.00){\makebox(0,0)[cc]{$_b$}}
\put(99.00,38.00){\makebox(0,0)[cc]{$_8$}}
\put(99.00,13.00){\makebox(0,0)[cc]{$_9$}}
\put(116.00,23.00){\makebox(0,0)[cc]{$_e$}}
\put(111.00,49.00){\makebox(0,0)[cc]{$_5$}}
\put(111.00,1.00){\makebox(0,0)[cc]{$_b$}}
\put(139.00,38.00){\makebox(0,0)[cc]{$_d$}}
\put(139.00,13.00){\makebox(0,0)[cc]{$_7$}}
\put(156.00,23.00){\makebox(0,0)[cc]{$_4$}}
\put(151.00,49.00){\makebox(0,0)[cc]{$_5$}}
\put(151.00,1.00){\makebox(0,0)[cc]{$_b$}}
\put(45.00,23.00){\makebox(0,0)[cc]{$_1$}}
\put(85.00,23.00){\makebox(0,0)[cc]{$_a$}}
\put(125.00,23.00){\makebox(0,0)[cc]{$_e$}}
\put(12.00,41.00){\makebox(0,0)[cc]{$_{F^d}$}}
\put(32.00,41.00){\makebox(0,0)[cc]{$_{D^3}$}}
\put(101.00,20.00){\makebox(0,0)[cc]{$_{B_0^2}$}}
\put(141.00,20.00){\makebox(0,0)[cc]{$_{A_1^6}$}}
\put(15.00,25.00){\makebox(0,0)[cc]{$_{D^4}$}}
\put(28.00,25.00){\makebox(0,0)[cc]{$_{F^9}$}}
\put(55.00,25.00){\makebox(0,0)[cc]{$_{C^1}$}}
\put(68.00,25.00){\makebox(0,0)[cc]{$_{E^b}$}}
\put(95.00,25.00){\makebox(0,0)[cc]{$_{F^9}$}}
\put(108.00,25.00){\makebox(0,0)[cc]{$_{D^2}$}}
\put(135.00,25.00){\makebox(0,0)[cc]{$_{E^b}$}}
\put(148.00,25.00){\makebox(0,0)[cc]{$_{C^1}$}}
\put(52.00,41.00){\makebox(0,0)[cc]{$_{C^3}$}}
\put(72.00,41.00){\makebox(0,0)[cc]{$_{A^0}$}}
\put(92.00,41.00){\makebox(0,0)[cc]{$_{B^e}$}}
\put(112.00,41.00){\makebox(0,0)[cc]{$_{C^2}$}}
\put(132.00,41.00){\makebox(0,0)[cc]{$_{A^5}$}}
\put(152.00,41.00){\makebox(0,0)[cc]{$_{C^d}$}}
\put(12.00,8.00){\makebox(0,0)[cc]{$_{A^5}$}}
\put(32.00,8.00){\makebox(0,0)[cc]{$_{E^f}$}}
\put(52.00,8.00){\makebox(0,0)[cc]{$_{A^4}$}}
\put(72.00,8.00){\makebox(0,0)[cc]{$_{D^4}$}}
\put(92.00,8.00){\makebox(0,0)[cc]{$_{F^a}$}}
\put(112.00,8.00){\makebox(0,0)[cc]{$_{F^5}$}}
\put(132.00,8.00){\makebox(0,0)[cc]{$_{F^0}$}}
\put(152.00,8.00){\makebox(0,0)[cc]{$_{F^1}$}}
\put(12.00,49.00){\makebox(0,0)[cc]{$_{A_0^2}$}}
\put(52.00,49.00){\makebox(0,0)[cc]{$_{A_0^6}$}}
\put(92.00,49.00){\makebox(0,0)[cc]{$_{A_0^f}$}}
\put(132.00,49.00){\makebox(0,0)[cc]{$_{A_0^c}$}}
\emline{177.00}{5.00}{155}{195.00}{5.00}{156}
\emline{197.00}{45.00}{157}{215.00}{45.00}{158}
\emline{177.00}{44.00}{159}{185.00}{36.00}{160}
\emline{187.00}{34.00}{161}{195.00}{26.00}{162}
\emline{195.00}{24.00}{163}{187.00}{16.00}{164}
\emline{187.00}{14.00}{165}{195.00}{6.00}{166}
\emline{186.00}{34.00}{167}{186.00}{16.00}{168}
\emline{197.00}{44.00}{169}{205.00}{36.00}{170}
\emline{205.00}{34.00}{171}{197.00}{26.00}{172}
\emline{197.00}{24.00}{173}{205.00}{16.00}{174}
\emline{176.00}{44.00}{175}{176.00}{6.00}{176}
\emline{206.00}{34.00}{177}{206.00}{16.00}{178}
\emline{207.00}{14.00}{179}{215.00}{6.00}{180}
\emline{216.00}{44.00}{181}{216.00}{6.00}{182}
\emline{176.00}{6.00}{183}{176.00}{5.00}{184}
\emline{176.00}{5.00}{185}{177.00}{5.00}{186}
\emline{215.00}{45.00}{187}{216.00}{45.00}{188}
\emline{216.00}{45.00}{189}{216.00}{44.00}{190}
\emline{185.00}{36.00}{191}{187.00}{34.00}{192}
\emline{205.00}{16.00}{193}{207.00}{14.00}{194}
\emline{195.00}{24.00}{195}{197.00}{26.00}{196}
\emline{205.00}{36.00}{197}{206.00}{35.00}{198}
\emline{206.00}{35.00}{199}{206.00}{34.00}{200}
\emline{187.00}{14.00}{201}{186.00}{15.00}{202}
\emline{186.00}{15.00}{203}{186.00}{16.00}{204}
\put(207.00,36.00){\vector(1,1){8.00}}
\put(187.00,36.00){\vector(1,1){8.00}}
\put(195.00,45.00){\vector(-1,0){18.00}}
\put(205.00,14.00){\vector(-1,-1){8.00}}
\put(185.00,14.00){\vector(-1,-1){8.00}}
\put(197.00,5.00){\vector(1,0){18.00}}
\put(177.00,48.00){\makebox(0,0)[cc]{$_{\beta_1}$}}
\put(197.00,48.00){\makebox(0,0)[cc]{$_{\beta_2}$}}
\put(196.00,2.00){\makebox(0,0)[cc]{$_{\beta_a}$}}
\put(183.00,35.00){\makebox(0,0)[cc]{$_{\beta_4}$}}
\put(202.00,35.00){\makebox(0,0)[cc]{$_{\beta_5}$}}
\put(183.00,15.00){\makebox(0,0)[cc]{$_{\beta_7}$}}
\put(202.00,15.00){\makebox(0,0)[cc]{$_{\beta_8}$}}
\put(216.00,48.00){\makebox(0,0)[cc]{$_{\beta_3}$}}
\put(216.00,2.00){\makebox(0,0)[cc]{$_{\beta_b}$}}
\put(196.00,33.00){\makebox(0,0)[cc]{$_{\alpha_1}$}}
\put(196.00,13.00){\makebox(0,0)[cc]{$_{\alpha_4}$}}
\put(211.00,23.00){\makebox(0,0)[cc]{$_{\alpha_3}$}}
\put(206.00,49.00){\makebox(0,0)[cc]{$_{\alpha_0}$}}
\put(206.00,1.00){\makebox(0,0)[cc]{$_{\alpha_5}$}}
\put(180.00,23.00){\makebox(0,0)[cc]{$_{\alpha_2}$}}
\put(196.00,21.00){\makebox(0,0)[cc]{$_{\beta_6}$}}
\put(187.00,49.00){\makebox(0,0)[cc]{$_{\beta_0}$}}
\put(177.00,2.00){\makebox(0,0)[cc]{$_{\beta_9}$}}
\end{picture}

\caption{Covering graph ${\Upsilon}_0$ of $e^b\cup d^2\cup c^1\cup
f^9-A_0$ and $\alpha$-$\beta$ denotations}
\end{figure}
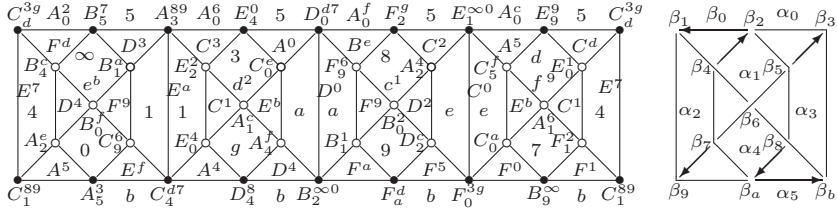

\noindent Each vertex of $\mathcal Y$ belongs to 12 copies of
$L(Q_3)$. Figure 7 shows the complements of vertex $A_0$ in 4 of the
12 copies of $L(Q_3)$ containing $A_0$, namely $e^b,d^2,c^1,f^9$,
which share the long vertical edges, successively present in the
copies $E^a,D^0,C^0,E^7$ of $K_4$, the last long vertical edge both
as the leftmost and rightmost edges in the shown covering graph, say
${\Upsilon}_0$, of $e^b\cup d^2\cup c^1\cup f^9-A_0$, where:

{\bf(a)} black vertices participate of the 8 $4$-holes containing
$A_0$, namely those labeled 5 on top and $b$ at the bottom; other
labels of 4-holes internal to them, respectively;

{\bf(b)} the labels $j$ of vertices $\Sigma_i$ appear as
superindices, as in $\Sigma_i^j$, (with $j$ also in the citations
$A_0^j$ of $A_0$ on top), or $\Sigma_i^{jj'}$, in case labels $j$
and $j'$ happen in contiguous copies of $L(Q_3)$;

{\bf(c)} each triangle contains the name $\Sigma^\ell$ of the copy
of $K_4$ containing it;

{\bf(d)} for each $\sigma^i=e^b,d^2,c^1,f^9$, the partition
$\Omega(\sigma^i)$ restricts as in the rightmost diagram, in which
darts indicate the first edges of monochromatic 2-paths whose final
vertex is $A_0$; as a result, the 4 mentioned long vertical edges
belong each to two different monochromatic 2-paths of contiguous
copies of $L(Q_3)$ in $\mathcal Y$;

{\bf(e)} alternate internal anti-diagonal monochromatic 2-paths
(i.e. from top-right to bottom-left) coincide with directions
reversed; (the middle vertices of these 4 2-paths are just two
neighbors of $A_0$ in $\mathcal S$, and their degree-1 vertices are
at distance 2 from $A_0$ in $\mathcal S$); and

{\bf(f)} the rightmost diagram contains denotations $\beta_i$,
($i\in[0,b]$), and $\alpha_j$, ($j\in[0,5]$), respectively for the
vertex and 4-hole labels in their positions in the 4 copies of
$L(Q_3)$.

\noindent Apart from the union $e^b\cup d^2\cup c^1\cup f^9$ of
copies of $L(Q_3)$ sharing $A_0$ in Figure 7, there are two other
unions of 4 copies of $L(Q_3)$ in $\mathcal Y$ sharing $A_0$. The
following display of the data in Figure 7 contains at its left the
$\alpha$-$\beta$ denotations of {\bf(f)}. Moreover, the data
corresponding to the 3 unions of 4 copies of $L(Q_3)$ sharing $A_0$
in $\mathcal Y$ are set (or encoded) in the arrays to the right and
below the $\alpha$-$\beta$ denotations (these solely for
$e^b,d^2,c^1,f^9$, respectively), where the leftmost array
summarizes ${\Upsilon}_0$, the two doubly repeated middle vertices
in ${\Upsilon}_0$ (as in {\bf(e)}) parenthesized to the right of
$A_0$ and the remaining data displayed in similar order, with the
two rightmost arrays preceded by the first one of their 4
corresponding $\alpha$-$\beta$ denotations, which condenses all
needed information of $\mathcal Y$ around $A_0$:
$$\begin{array}{llllll}
^{\alpha_0\beta_0\alpha_5=52b,56b,5fb,5cb}
_{\beta_1\beta_2\beta_3=g78,90d,7g\infty,093}\!&\!
^{A_0(B_0A_1)}_{(E^7e^bE^ad^2D^0c^1C^0f^9)} \!&\!^{f62}_{41g}\!&\!
^{A_0(A_1A_g)}_{(E^7d^gC^0d^1E^ae^9D^0e^8)}\\
^{\beta_4\alpha_1\beta_5=c\infty a,23e,684,fd1}
_{\alpha_2\beta_6\alpha_3=4f1,1ca,a2e,e64}\!&\!
^{(C_dB_5A_3E_4D_0F_2E_1E_9)}_{(B_4B_1E_2C_0F_9A_2C_5E_0)}
\!&\!^{5d9}_{8c0}\!&\!
^{(C_1E_eE_gD_3A_3C_7B_fC_6)}_{(D_1A_fC_3B_1C_2B_gB_eA_2)}\\
%\vspace*{1mm}
^{\beta_7\alpha_4\beta_8=e06,4gf,19c,a72}
_{\beta_9\beta_a\beta_b=93d,78\infty,0d3,g\infty 8}\!&\!
^{(A_2C_9E_0A_4B_1D_2C_0F_1)}_{(C_1A_5C_4D_4B_2F_aF_0B_9)}
\!&\!^{\infty 7b}_{3ae}\!&\!
^{(B_gC_eA_2D_gA_fB_3B_1C_f)}_{(A_eD_eE_1E_3C_gC_bB_2C_a)}\\
\!&\!\!&\!\!&\!\!&\!\!&\!\!\\
\!&\!\!&\!^{cb6}_{804}\!&\!
^{A_0(A_gB_0)}_{(E^7d^fD^0c^gC^0f^8E^ae^6)} \\
\!&\!\!&\!^{fe3}_{g57}\!&\!
^{(A_eE_dD_0F_fE_gE_8C_4B_c)}_{(E_fC_0F_8A_fC_cE_0B_dB_g)} \\
\!&\!\!&\!^{d12}_{a9\infty}\!&\!
^{(E_0A_dB_gD_fC_0F_gA_fC_8)}_{(C_dD_dB_fF_7F_0B_8C_gA_c)} \\
\end{array}$$
\noindent Some edges are shared by two of these 3 unions. In fact,
each of the edges bordering the central 2-paths $\omega$ in
anti-diagonal 4-paths in ${\Upsilon}_0$ is present also in one of
the two covering graphs, say ${\Upsilon}_1$ and ${\Upsilon}_2$,
corresponding to the two rightmost arrangements above, one encoded
on top and the other at the bottom of the display, respectively. For
example, the edge $B_1A_3$ of $e^b$ on ${\Upsilon}_0$ appears in
${\Upsilon_1}$. Also, the labels
$\{\alpha_0\alpha_4,\alpha_1\alpha_5,\alpha_2\alpha_3\}$ of opposite
copies of $L(Q_3)$, just sharing vertex $A_0$, are images of
vertices at distance 3 in $\mathcal S$ via the duality $\phi$ (but
copies of $L(Q_3)$ sharing a triangle containing $A_0$ are images of
vertices at distance 7). The following permutations on the set
$\{\alpha_0,\ldots,\alpha_5,\beta_0,\ldots,\beta_{11}\}$ relate the
labels of the 12 copies of $L(Q_3)$ sharing $A_0$:
$$^{e^b\rightarrow d^2\rightarrow c^1\rightarrow f^9\rightarrow
e^b:}_{\hspace*{26mm}(\alpha_0)(\alpha_5)(\beta_0 \beta_4 \beta_6
\beta_8)(\beta_1 \alpha_4 \beta_2 \beta_9)(\beta_3 \beta_a \alpha_1
\beta_b)(\beta_5 \alpha_3 \alpha_2 \beta_7);}$$ $$^{
e^bd^2c^1f^9\rightarrow d^gd^1e^9e^8\rightarrow
d^fc^gf^8e^6\rightarrow e^bd^2c^1f_9:}_{\hspace*{26mm}(\alpha_0
\beta_4 \beta_6)(\beta_0 \alpha_5 \beta_8)(\beta_1 \beta_3
\alpha_2)(\beta_2 \alpha_4 \alpha_3)(\alpha_1 \beta_7
\beta_b)(\beta_5 \beta_a \beta_9).}$$ \noindent The following
permutations allow to relate the labels of the 12 cuboctahedral
subgraphs sharing $A_0$ to those sharing $B_0,C_0,D_0,E_0,F_0$:
$$\begin{array}{l}
^{A_0\rightarrow
B_0\;:\;(\alpha_0\alpha_3\beta_a\alpha_1\beta_5\alpha_4\beta_7\beta_b\beta_2\beta_1\beta_3\beta_0\alpha_5\beta_4\beta_8\alpha_2\beta_9\beta_6);}
_{A_0\rightarrow C_0\;:\;(\alpha_0\beta_1\beta_2\beta_0\alpha_4\beta_3)(a1\beta_9\beta_6)(\alpha_2\beta_a\beta_4)(a3\beta_7\alpha_5)(\beta_5\beta_8\beta_b);}\\
^{A_0\rightarrow
D_0\;:\;(\alpha_0\beta_8\alpha_2\beta_0\beta_a\beta_b\beta_6\beta_5\beta_4)(\alpha_1\beta_9\beta_7a3\alpha_4\beta_3\beta_2\alpha_5\beta_1);}
_{A_0\rightarrow E_0\;:\;(\alpha_0\beta_b\beta_0\beta_a\beta_8\alpha_2\beta_6\beta_3\beta_2\alpha_4\beta_5\beta_1)(\alpha_1\beta_7\alpha_3\alpha_5\beta_4)(\beta_9);}\\
^{A_0\rightarrow
F_0\;:\;(\alpha_0\beta_b\alpha_4\beta_3\beta_5\alpha_2\alpha_1\beta_9\beta_a)(\alpha_3\beta_0\beta_2\beta_1\beta_6\beta_7\beta_8\alpha_5\beta_4).}
\end{array}$$
\noindent Additions mod 17 yield the remaining information for
copies of $K_4$ and $L(Q_3)$ neighboring each vertex of $\mathcal
Y$. In sum, we have the following theorem.%$\mathcal Y$ has the following properties.

\begin{thm}
In addition to {\rm Theorem 4.1}, the following properties of $\mathcal Y$ hold:

\noindent{\bf(1)} $\mathcal Y$ is a connected union of $102$ copies
$\sigma$ of $L(Q_3)$, each with an edge partition $\Omega(\sigma)$
into $2$-paths;

\noindent{\bf(2)} each edge in $\mathcal Y$ is shared exactly by $4$
copies of $L(Q_3)$ in $\mathcal Y$;

\noindent{\bf(3)} each copy $\Delta$ of $K_3$ $($resp. each $2$-path
$\omega\in \Omega(\sigma))$ in a copy $\sigma$ of $L(Q_3)$ in
$\mathcal Y$ is shared exactly by two copies $\sigma,\sigma'$ of
$L(Q_3)$ in $\mathcal Y$;

\noindent{\bf(4)} Each two copies of $L(Q_3)$ sharing a copy
$\Delta$ of $K_3$ in ${\mathcal Y}$ share $\Delta$ with exactly one
copy of $K_4$ in $\mathcal Y$;

\noindent{\bf(5)} each $4$-hole in $\mathcal Y$ happens in just one
copy of $L(Q_3)$ in $\mathcal Y$;

\noindent{\bf(6)} $\mathcal Y$ is an $\Omega$-preserving
$\{L(Q_3)\}_{K_3}$-{\rm UH} graph;

\noindent{\bf(7)} $\mathcal Y$ is $\{K_4,L(Q_3)\}_{K_3}$-{\rm UH};

\noindent{\bf(8)} the vertices and copies of $L(Q_3)$ in $\mathcal
Y$ are the points and lines of a self-dual $1$-configuration
$(102_4)_1$ whose Menger graph is again $\mathcal Y$.
\end{thm}

\noindent In  Theorem 6.1{\bf(3)}, for each triangle $\Delta$ in
$\sigma$, the copies $\sigma,\sigma'$ of $L(Q_3)$ intersect exactly
in $\Delta$, while for each $2$-path $\omega\in\Omega(\sigma)$ in
$\sigma$, not only $\omega$ is shared by $\sigma,\sigma'$, but these
also share a vertex at distance 2 from the ends of $\omega$. This
common distance, 2, is realized by 2-paths in the other two colors
distinct from the color of $\omega$, in each of $\sigma$ and
$\sigma'$, as in Figure 4, where for example the dark-gray-colored
2-path $F_4D_2B_4$ (present both in $a^0$ and $c^3$) is at distance
2 from vertex $D_4$ (also present in $a^0$ and $c^3$) via the
black-colored path $B_4F_dD_4$ and the light-gray-colored path
$F_4C_4D_4$.

\proof It only remains to prove item {\bf(8)}. We explain how a
monochromatic 2-path-preserving isomorphism
$\Psi':\sigma'_1\rightarrow\sigma'_2$ between two copies of $L(Q_3)$
$\sigma'_1,\sigma'_2$ in ${\mathcal Y}$ extends to an automorphism
of $\mathcal S$. Both $\sigma'_1$ and $\sigma'_2$ are colored as in
Figure 4 with $\Psi'$ respecting the color structure, thus inducing
a 1-1 correspondence between the color classes of $\sigma'_1$ and
$\sigma'_2$. In each copy of $L(Q_3)$ in $\mathcal Y$ there are
exactly 12 monochromatic 2-paths, 4 in each of the 3 colors, and
exactly 12 dichromatic 2-paths not contained in any triangle, a
total of 24 2-paths not contained in any triangle. A
$\Psi':\sigma'_1\rightarrow\sigma'_2$ as mentioned can be extended
to an automorphism of $\mathcal Y$ because the information gathered
in $\sigma'_i$ comes via sextets from corresponding information in a
subgraph $\overline{\sigma'}_i$ of $\mathcal S$, ($i=1,2$), so that
$\Psi'$ arises from an isomorphism
$\overline{\Psi'}:\overline{\sigma'}_1\rightarrow\overline{\sigma'}_2$.
However, $\overline{\sigma'}_i=\overline{\sigma}_i$, ($i=1,2$), for
a corresponding copy $\sigma_i$ of $L(Q_3)$ in $\mathcal Y$, but
while the vertices of $\sigma'_i$ are denoted like the degree-2
vertices of $\overline{\sigma'}_i=\overline{\sigma}_i$, the vertices
of $\sigma_i$ are denoted like the degree-3 vertices of
$\overline{\sigma}_i=\overline{\sigma'}_i$. Here the pairs
$(\sigma_i,\sigma'_i)$ are of the form $(\Sigma^j,\sigma^j)$, where
$(\Sigma,\sigma)\in\{(A,a),(B,b),(C,c),(D,d),(E,e),(F,f)\}$ and
$j\in\Z_{17}$. Then
$\overline{\Psi'}=\overline{\Psi}:\sigma_1\rightarrow\sigma_2$ is a
corresponding map as in the proof of Theorem 4.1. But now
$\overline{\Psi'}=\overline{\Psi}$ extends to an automorphism of
$\mathcal S$. This takes us to an automorphism of $\mathcal Y$ that
extends $\Psi'$, as claimed above.

\noindent For example, the black 2-path $B_4F_dD_4$ in the copy
$a^0$ of $L(Q_3)$ in $\mathcal Y$ represented in Figure 4 arise from
the 3-paths $B_4E_4F_4F_d$ and $F_dF_4E_4D_4$ in $\mathcal S$, which
share the 2-path $F_dF_4E_4$ and differ otherwise, so their union
$(B_4E_4F_4F_d)\cup(F_dF_4E_4D_4)$ is realized by a tree $T_1$ with
just one vertex of degree 3, namely $E_4$, from which two 1-paths
and one 2-path depart. A similar tree $T_2$ is obtained from the
black 2-path $D_dF_4B_d$ in Figure 4. But $T_1\bigcap T_2=F_dF_4$, a
terminal 1-path of $T_i$ on its 2-path departing from $t_i$, for
both $i=1,2$, where $t_1=E_4$ and $t_2=E_d$, the vertex of degree 3
in $T_2$. The other two black 2-paths in Figure 4 behave similarly,
leading to trees $T_3$ and $T_4$ intersecting at the 1-path
$B_0E_0$. Similar behavior holds for the dark gray and the light
gray quadruples of 2-paths in Figure 4, leading to pairs of trees
that intersect respectively at the 1-paths $D_4D_2$, $B_dC_d$ and
the 1-paths $B_4C_4$, $D_fD_d$. Thus, if $\sigma'_1$ is this copy of
$L(Q_3)$ in $\mathcal Y$, then $\overline{\sigma'}_1$ coincides with
$\overline{\sigma}_1$, where $\sigma_1=A^0$.\qfd

\section{Using the Biggs-Smith association scheme}

\noindent The 2-paths $\omega$ of Theorem 6.1{\bf(3)} rearrange into
an edge partition $\mathcal I$ of $\mathcal Y$ into 102\, $4$-holes.
In fact, each 4-hole in $\mathcal I$ is the union of 4 successive
2-paths $\omega_0,\omega_1,\omega_2,\omega_3$ from 4 respective
partitions
$\Omega(\sigma^0),\Omega(\sigma^1),\Omega(\sigma^2),\Omega(\sigma^3)$
of $L(Q_3)$ into 2-paths, with each two successive 2-paths
$\omega_i,\omega_{i+1}$ here overlapping in just one edge, (subindex
addition taken mod 4).

\noindent $\mathcal I$ can be reconstructed by adding $r\in\Z_{17}$
uniformly mod 17 to all indexes in the following generating-set
table of its member 4-holes, from those 4-holes shown in the left
column of the table. In each line of the table, the 4 pairs of
copies $\sigma^i_j$ of the disconnected graph $4P_3$ shown to the
right (as in (7) above) overlap at succeeding pairs of 2-paths of
the 4-hole shown on their left. This is continued to its right by
the citation of two vertices that alternatively are at distance 2
from the ends of those composing 2-paths:
$$\begin{array}{||l||c|c|c|c||}\hline
^{(A_2B_0B_1A_g)\,\,A_0A_1}_{(C_0A_gE_0A_1)\,\,A_0B_0}&^{(c^1_3\,\,e^b_2)}_{(d^f_2\,\,f^8_1)}&^{(e^7_2\,\,c^0_2)}_{(c^0_1\,\,d^0_1)}&^{(d^1_3\,\,e^8_3)}_{(d^2_3\,\,f^9_1)}&^{(e^a_3\,\,d^0_2)}_{(e^7_1\,\,e^a_1)}\\
^{(C_4E_0C_dA_0)\,\,B_0C_0}_{(D_0A_0F_0C_0)\,\,B_0E_0}&^{(a^0_1\,\,f^0_1)}_{(c^g_2\,\,c^1_3)}&^{(f^9_2\,\,d^f_1)}_{(f^8_2\,\,f^9_3)}&^{(e^6_2\,\,e^b_2)}_{(a^d_2\,\,a^4_3)}&^{(d^2_1\,\,f^8_3)}_{(d^f_3\,\,d^2_2)}\\
^{(C_8B_0B_4C_d)\,\,\,C_0C_4}_{(D_4D_fE_2E_0)\,D_0D_2}&^{(a^4_3\,\,e^a_1)}_{(a^2_2\,\,b^g_3)}&^{(e^b_1\,\,a^0_2)}_{(b^e_2\,\,d^0_3)}&^{(f^4_3\,\,e^f_3)}_{(d^2_2\,\,b^5_2)}&^{(e^6_3\,\,f^0_2)}_{(b^3_3\,\,a^0_3)}\\
^{(F_0D_2B_0D_f)\,\,D_0E_0}_{(F_8B_0F_9D_0)\,\,\,E_0F_0}&^{(c^1_1\,\,a^4_2)}_{(c^0_1\,\,f^0_1)}&^{(a^0_3\,\,d^0_1)}_{(c^1_2\,\,a^d_1)}&^{(a^d_3\,\,c^g_1)}_{(b^5_2\,\,b^c_2)}&^{(b^3_1\,\,b^e_1)}_{(a^4_1\,\,c^g_3)}\vspace*{0.5mm}\\
^{(E_8\,E_0F_gF_9)\,\,F_0F_8}&^{({b^3_1}\,\,{f^8_2})}&^{({b^c_3}\,\,{c^0_2})}&^{({c^8_3}\,\,{b^d_3})}&^{({f^0_3}\,\,{b^5_1})}\\\hline
\end{array}$$
The vertices of each such 4-hole coincide in notation with the
degree-1 vertices of a tree $T$ in $\mathcal S$ isomorphic to
$T_0^\infty$, (itself present in the 4th row of this table), with
the two vertices that follow each 4-hole being the vertices of
degree 3 in $T$. These data insure that $\mathcal Y$ is
${\mathcal I}$-UH.

\noindent Of the 24 2-paths in a copy $\sigma^i$ of $L(Q_3)$ in
$\mathcal Y$, 12 are in the partition $\Omega(\sigma^i)$ of
$\sigma^i$. The other 12 form a different edge partition
$\Omega'(\sigma^i)\ne \Omega(\sigma^i)$ of $\sigma^i$. The
family of 2-paths in all of the $\Omega'(\sigma^i)$\thinspace s
reassembles, by means of unions of those of its members having a
common degree-2 vertex, as a family $\mathcal J$ of 306 copies of
$K_{1,4}$.

\noindent A generating-set table for $\mathcal J$ representing 18
copies of $K_{1,4}$ is shown subsequently, with the remaining copies
of $K_{1,4}$ obtained from those 18 by uniform addition of
$r\in\Z_{17}$ to all indexes $i\in\Z_{17}$ of vertices $\Sigma_i$
and subgraphs $\sigma^i_j$, where $j=1,2,3$ stands for black, dark
gray and light gray, respectively. This generating-set table has
each entry starting with a vertex $\Sigma_0$ of degree 4 in a copy
of $K_{1,4}$ in $\mathcal J$ followed by 4 parenthesized
expressions, each containing as its central entry a neighbor
$\Sigma'$ of $\Sigma_0$ flanked by two subgraphs $\sigma^i_j$ to
which the edge $\Sigma_0\Sigma'$ belongs, so that each participating
$\sigma^i$ appears repeated twice --- with 2 different colors
$j,j'$, as $\sigma^i_j$ and $\sigma^i_{j'}$
--- once before a right parenthesis and once after the subsequent
left parenthesis, the first of the 4 left parentheses considered
subsequent to the last right parenthesis, in a mod 4 fashion:
$$\begin{array}{|l|}\hline
^{A_0\,(e^b_3 \, A_3 \, d^1_2)\,\,(d^1_1 \, E_1 \, c^1_1)\,(c^1_2 \,B_2 \, e^8_3)\,\,(e^8_1 \, C_1 \, e^b_1)}
_{A_0\,(f^8_3 \, C_4 \,d^2_1)\,(d^2_2 \, D_0 \, d^f_3)(d^f_1 \, C_d \, f^9_2)\,(f^9_3 \,F_0 f^8_2)}\\
^{A_0\,(d^g_3 \, A_e \, e^6_3)\,\,(e^6_1 \, C_g \, e^9_1)\,(e^9_2 \,B_f \, c^g_3)\,\,(c^g_2 \, E_g \, d^g_1)}
_{B_0\,(e^6_1 \, B_d \,a^0_3)\,\,(a^0_2 \, B_4 \, e^b_1)\,(e^b_3 \, C_9 \, f^0_3)\,\,(f^0_2\, C_8 \, e^6_3)}\\
^{B_0\,(e^7_3 \, A_f \, d^0_3)\,\,(d^0_2 \, A_2 \, e^a_3)\,(e^a_2 \,B_g \, c^0_3)\,\,(c^0_2 \, B_1 \, e^7_2)}
_{B_0\,(a^4_2 \, D_2 \,c^1_1)\,\,(c^1_2 \, F_9 \, a^d_1)\,(a^d_3 \, D_f \, c^g_1)\,\,(c^g_3\, F_8 \, a^4_1)}\\
^{C_0\,(d^f_3 \, D_0 \, d^2_2)\,\,(d^2_3 \, A_1  f^9_1)\,(f^9_3 \,F_0 \, f^8_2)\,\,(f^8_1  A_g  d^f_2)}
_{C_0\,(e^7_2 \, A_d \,e^2_2)\,\,\,(e^2_1 \, B_9 \, a^d_3)\,(a^d_1 \, E_d \,f^d_1)\,\,(f^d_3 \, C_5  e^7_3)}\\
^{C_0\,(d^f_1 \, B_8 \, a^4_2)\,\,(a^4_1 \, E_4 \, f^4_1)\,(f^4_2 \,C_c \, e^a_3)\,\,(e^a_2 \, A_4  d^f_2)}
_{D_0\,(b^c_1 \, F_f \, b^1_1)\,\,\,(b^1_2 \, E_d \,d^f_2)\,\,(d^f_1  B_f \, a^f_1)\,(a^f_2 \, D_b \, b^c_3)}\\
^{D_0\,(a^d_1 \, F_9 \, c^1_2)\,\,\,(c^1_3 \, A_0 \, c^g_2)\,\,(c^g_3 \, F_8 \, a^4_1)\,\,\,(a^4_3 \, C_0 \, a^d_2)}
_{D_0\,(b^5_3 \, D_6 \, a^2_3)\,\,(a^2_1 \, B_2\, d^2_1)\,\,(d^2_3 \, E_4 \, b^g_2)\,\,(b^g_1 \, F_2 \, b^5_1)}\\
^{E_0\,(a^0_2 \, D_d \, b^e_2)\,\,\,(b^e_2 \, E_2 \, d^0_3)\,\,(d^0_2 \, E_f \, b^3_2)\,\,(b^3_3 \, D_4 \, a^0_3)}
_{E_0\,(b^5_3 \, F_1 \, c^0_3)\,\,\,\,(c^0_2 \,F_g \, b^c_3)\,\,\,(b^c_1 \, E_9 \, f^0_2)\,\,\,(f^0_3 \, E_8 \,b^5_1)}\\
^{E_0\,(f^9_1 \, A_1  d^2_3)\,\,\,(d^2_1  C_4 \, f^8_3)\,\,\,(f^8_1 \, A_g  d^f_2)\,\,(d^f_1 \, C_d \, f^9_2)}
_{F_0\,(c^g_2 \, A_0 \,c^1_3)\,\,\,(c^1_1 \, D_2 \, a^4_2)\,\,\,(a^4_3 \, C_0 \,a^d_2)\,\,(a^d_3 \, D_f \, c^g_1)}\\
^{F_0\,(b^d_2 \, D_8 \, b^3_2)\,\,\,(b^3_3 \, F_7 \, c^8_2)\,\,\,\,(c^8_1 \, B_8 \, f^8_1)\,\,\,(f^8_3 \, E_g \, b^d_1)}
_{F_0\,(f^9_2 \, E_1 \,b^4_1)\,\,\,(b^4_2 \, D_9 \, b^e_2)\,\,\,(b^e_3 \, F_a \,c^9_3)\,\,\,\,(c^9_1 \, B_9 \, f^9_1)}\\\hline
\end{array}$$
Here, a copy of $K_{1,4}$ with degree-4 vertex $\Sigma_i$ has its
degree-1 vertices as those of a binary tree of $\mathcal S$ with
depth 2 and whose root is one of the 3 neighbors of $\Sigma_i$.
Thus, there are 3 such copies of $K_{1,4}$. As a result, in contrast
to the fact mentioned above that $\mathcal Y$ is $\mathcal I$-UH,
now any homomorphism between members of $\mathcal J$ preserving the
order of presentation of the degree-1 vertices in corresponding
copies of $K_{1,4}$, as in the table above (with the expressed
parenthetical behavior with respect to the $\sigma_j^i$\thinspace
s), extends to an automorphism of $\mathcal Y$. On the other hand,
each copy $\sigma$ of $L(Q_3)$ in $\mathcal Y$ intersects 8 other
copies of $L(Q_3)$ in a triangle each, and 12 other copies of
$L(Q_3)$, each in a 2-path of $\Omega(\sigma)$ and one more vertex
at distance $2$\, from the ends of the 2-path.

\noindent The graph ${\mathcal I}\,'$ generated by the (diagonal)
chords of the 4-cycles of $\mathcal I$ coincides with ${\mathcal
S}_2$. On the other hand, by expressing the copies of $K_{1,4}$ in
$\mathcal J$ as $u(v)(w)(x)(y)$, (for example the copy of $K_4$ in
the first line of the last table as $A_0(A_3)(E_1)(B_2)(C_1)$), we
consider the graph ${\mathcal J}\,'$ generated by the corresponding
4-cycles $(v,w,x,y)$. Then ${\mathcal J}\,'$ coincides with
${\mathcal S}_4$. We obtain the following final result.

\begin{thm} Using the Biggs-Smith association scheme, it is obtained that
 ${\mathcal Y}={\mathcal S}_3$.\end{thm}

\proof As ${\mathcal I}\,'={\mathcal S}_2$ and ${\mathcal
J}\,'={\mathcal S}_4$, and because $\mathcal S$ has girth 9 and
$\mathcal Y$ was constructed from the family $({\mathcal C}_9)_3$ of
distance-3 digraphs of directed $9$-cycles in the set ${\mathcal
C}_9$ of 136 directed 9-cycles in Section 3, taking into account the
discussion previous to the statement, we arrive at
$$K_{102}\,=\,{\mathcal S}\cup{\mathcal S}_2\cup
{\mathcal S}_3\cup{\mathcal S}_4\,=\,{\mathcal S}\cup{\mathcal
I}'\cup {\mathcal Y}\cup{\mathcal J}',$$ and so ${\mathcal
Y}={\mathcal S}_3$. \qfd

\end{document}